%% file: main.tex
\documentclass[12pt]{article}
\include{macros}

\usepackage{mathrsfs}
\usepackage{geometry}  
\usepackage{amssymb}
\usepackage{amsmath}
\usepackage{xcolor}
\usepackage{subcaption}
\usepackage[ruled,vlined]{algorithm2e}
\usepackage{graphicx}
\usepackage{soul}

\newcommand{\LL}[1]{\textcolor{black}{#1}} 

\begin{document}
%
\title{Numerical wave propagation aided by deep learning}


\author{Hieu Nguyen\footnote{Department of Mathematics and Computer Science, University of Basel, Switzerland}\quad and\quad Richard Tsai\footnote{ Department of Mathematics and Oden Institute for Computational Engineering and Sciences, The University of Texas at Austin, Austin TX, USA}}

\date{}
\maketitle
\begin{abstract}
We propose a deep learning approach for wave propagation in media with multiscale wave speed, using a second-order linear wave equation model.
We use neural networks to enhance the accuracy of a given inaccurate coarse solver, which under-resolves a class of multiscale wave media and wave fields of interest. Our approach involves generating training data by the given computationally efficient coarse solver and another sufficiently accurate solver, applied to a class of wave media (described by their wave speed profiles) and initial wave fields. We find that the trained neural networks can approximate the nonlinear dependence of the propagation on the wave speed as long as the causality is appropriately sampled in training data. 
We combine the neural-network-enhanced coarse solver with the parareal algorithm and demonstrate that the coupled approach improves the stability of parareal algorithms for wave propagation and improves the accuracy of the enhanced coarse solvers.
\end{abstract}

\section{Introduction}

The computation of wave propagation is a critical ingredient in many important applications.
For example, in inverse problems involving the wave equation as the forward model, efficient numerical simulation of wave propagation directly corresponds to a more precise and faster imaging pipeline.

In this work, we will focus on the model initial value problem:
\begin{align} \label{eq:waveequation} 
    &u_{tt} = c^2(x)\Delta u,\,\,\, x\in[-1,1)^2,0\le t <T, \\ 
    &u(x,0) =u_0(x), \nonumber\\
    &u_t(x,0) =p_0(x), \nonumber
\end{align} 
with periodic boundary conditions. 
Furthermore, we shall assume that (i) the wave speed $c(x)$ is only piecewise smooth with non-trivial total variation, and (ii) the wavelengths in the initial wave field are at least one order of magnitude smaller than the domain size. 

Performing large-scale high-fidelity simulations of wave propagation in heterogeneous media is a challenging task. Standard numerical methods typically require very fine spatio-temporal grids to fully resolve the smallest wavelength in the solutions, the geometry in the wave speed's discontinuities, and to reduce excessive numerical dispersion errors. Furthermore, the stability and accuracy requirement imposes an additional restriction on the time step size, making long-time simulation even more difficult. 

Many existing multiscale methods tackle the wave equation by assuming a separation of scales in the medium's wave speed. The assumption on more global structures such as periodicity or stationary ergodicity is also common.
Such assumptions allow a cleverly designed algorithm to achieve a computational complexity that is essentially independent of the need to resolve the smallest scale globally in space-time.
These are algorithms may employ either online or offline computations to approximate the effective wave propagation accurately. We mention the Heterogeneous Multiscale Methods (HMM), see e.g.
\cite{weinan2003heterognous}\cite{abdulle2011finite}\cite{engquist2012multiscale}, which typically solve the effective problem, free of small parameters, defined by localized online computations. 
The Multiscale Finite Element Method (MsFEM) style algorithms \cite{efendiev2009multiscale}, including the Localized Orthogonal Decomposition (LOD) method, e.g.  \cite{abdulle2017localized}, compute specialized basis functions for each given multiscale medium. See also \cite{fu2019high}, and \cite{owhadi2008numerical}. These methods do not require scale separation, but need to fully resolve the medium by the basis functions.
As the basis functions are medium-specific, these algorithms typically aim to solve the same equation repeatedly with different initial conditions. See also \cite{fu2019high}, and \cite{owhadi2008numerical}.
It is also possible to construct low-rank approximations of the numerical wave propagator using data.
Druskin et. al. \cite{druskin2018nonlinear,borcea2020reduced} introduced a data-driven reduced order model (ROM) for the wave equation. Specifically, using wave field time snapshots the authors construct a wave propagator allowing very accurate inversion of the medium. 

\LL{ In this paper, we will tap into enhancing the computation on a coarse grid by a deep learning approach and combine it with a time-parallelization algorithm, called parareal, to reduce the overall computational time on the wall-clock. 
The parareal scheme \cite{parareal-LMT01} relies on a efficient coarse solver and an accurate but expensive fine solver. It is conceived with the target to be scalable on high performance computing (HPC) clusters. The fine solver is to be computed for shorter time intervals and in parallel.
This type of schemes introduce additional iterations to gain concurrency on the temporal dimension, making it a natural combination with the classical domain decomposition technique to fully utilize the computing power. There have been several works on efficient implementations of spatial-temporal domain decomposition \cite{d2022scalable}\cite{emmett2014efficient}\cite{minion2015interweaving}. Their results on the strong and weak scaling of the scheme demonstrated the scalability on HPC. Some simulation examples in three-dimension are \cite{croce2014parallel},\cite{speck2014space}. 
A scalable open-source package MGRIT (multigrid in time) has been developed by Falgout \textit{et. al.} \cite{falgout2014parallel}.}

\LL{A key limitation of parareal scheme is the convergence property for convection-dominated problems; see e.g. \cite{ANT-thetaparareal}\cite{mikioInfluence18}\cite{Ruprecht18}. The underlying reason is the lack of dissipation and
the numerical dispersion errors produced by the coarse solver. This instability issue has prevented the adoption of the scheme in solving large-scale hyperbolic problems. One approach to stabilize the error convergence is proposed in \cite{nguyen2020stable}, in which the computation of a coarse, under-resolving solver is corrected for phases using data computed ``onine'' and in-parallel by an accurate solver. In this manuscript, we also develop a stabilized parareal scheme, but we use deep convolution neural network trained ``offline'' with a predetermined set of data.
By the nature of convolution, a trained network can be used in computational involving larger grids. 
}
\LL{\paragraph{Our contributions.} This paper serves as a proof-of-concept of using machine learning approaches to enhance coarse computations by classical numerical methods.
Our main contributions include two areas:\\ 
(1) Introducing a working deep learning approach in a supervised learning setup that enhances the accuracy of a lower fidelity computation performed on a coarse grid. In this regard, the contribution includes the introduction of a novel wave energy based optimization model and a systematic way to generate training data that properly samples the strong causality in the wave propagation. The resulting neural networks significantly reduce the numerical dispersion error and the loss of high frequency modes and energy in the coarse computations.\\
(2) Enabling stable parallel-in-time computations for wave propagation, which further improves the accuracy of solutions enhanced by the deep learning approach. }


\subsection{The proposed approach}
Our approach departs from the above multiscale methodologies.  
We propose a supervised deep learning framework to enhance the accuracy of a low-fidelity and computationally cheaper solver for~\eqref{eq:waveequation}. 
After describing how neural networks are used in our approach, we will present our choice of neural network model (the architecture, the input, and output variables), and our approach to generating training data
that consistently samples the causality in wave propagation in a class of media with piecewise smooth wave speeds.
We will also discuss the rationales behind our approach and report several representative numerical studies.

We assume  a time-reversible coarse solver, $\mathcal{G}_{\Delta t}\mathbf{u}\equiv\mathcal{G}_{\Delta t}[\mathbf{u}, c]$ and 
an accurate fine solver, $\mathcal{F}_{\Delta t}\mathbf{u}\equiv\mathcal{F}_{\Delta t}[\mathbf{u}, c]$. 
These solvers advance a wave field 
$\mathbf{u}(x,t)=(u, \partial_t u)$ to time $t+\Delta t,$ on a medium with the given wave speed profile $c(x)$.
 \LL{We presume these explicit solvers are stable. The coarse solver is computationally cheap, but may under resolve the wave fields and the wave speed.}

The fine solver is used to generate data for improving the coarse solver. \LL{There is essentially no restriction to the choice of fine solvers, except that it has to be stable and sufficiently accurate for propagation in the class of wave speeds under consideration.}
As these two solvers operate on grids of different resolutions (a coarse grid and a fine grid), we will extend the coarse solver to one defined on the same grid as the fine solver for comparing the wave fields propagated by the two solvers. The extension is defined via a restriction (projection) operator $\mathcal{R}$, which maps fine grid functions to the coarse grid. We shall also need a prolongation (e.g. interpolation) operator $\mathcal{I}$, which maps coarse grid functions to the fine grid. 

We reorganize the time stepping by the fine solver into
\begin{equation}
\begin{aligned}
    \mathbf{u}_{n+1}:= \mathcal{F}_{\Delta t} \mathbf{u}_n&= \mathcal{F}_{\Delta t}\left(\mathcal{IR} \mathbf{u}_n\right)  + \left[\mathcal{F}_{\Delta t} \mathbf{u}_n - \mathcal{F}_{\Delta t}\left( \mathcal{IR} \mathbf{u}_n\right) \right].
\end{aligned}
\end{equation}
In general, $\mathcal{R}$ acts as a low-pass filter, and 
$\mathcal{IR}\mathbf{u} \neq \mathbf{u}$. 
Furthermore, it is reasonable to approximate 
$\mathcal{F}_{\Delta t}\left(\mathcal{IR} \mathbf{u}_n\right)$ by a computationally cheaper strategy using a coarser grid --- we consider $\mathcal{I}\mathcal{G}_{\Delta t}\left(\mathcal{R}\mathbf{u}_n\right)$.
The second term above, $\mathcal{F}_{\Delta t} \mathbf{u}_n - \mathcal{F}_{\Delta t}\left( \mathcal{IR} \mathbf{u}_n\right)$, typically consists mostly the higher frequency components of $\mathcal{F}_{\Delta t}\mathbf{u}.$


If we further introduce a fixed-point iteration that allows the fine solver to be applied in parallel, we would derive a parareal scheme \cite{parareal-LMT01}:
\begin{equation}\label{eq:vanilla-parareal-itrs}
\begin{aligned}
    \mathbf{u}_{n+1}^{k+1} &:=  \mathcal{I}\mathcal{G}_{\Delta t} \mathcal{R} \mathbf{u}_n^{k+1} + \left[\mathcal{F}_{\Delta t} \mathbf{u}_n^{k} - \mathcal{I}\mathcal{G}_{\Delta t} \mathcal{R} \mathbf{u}_n^{k}\right],~~~k=0,1,\dots,K,
\end{aligned}
\end{equation}
with
\begin{equation}
    \mathbf{u}_{n+1}^0:=\mathcal{I}\mathcal{G}_{\Delta t} \mathcal{R}\mathbf{u}_n^0,~~~n=0,1,2,\dots.
\end{equation}
The second term on the right hand side of \eqref{eq:vanilla-parareal-itrs}
serves as an attempt to add back the missing high frequency components, using information from the previous iteration.

The iteration in \eqref{eq:vanilla-parareal-itrs} is often unstable when the solvers have too little dissipation.
In such cases, the stability and convergence of \eqref{eq:vanilla-parareal-itrs}
will rely heavily on $\mathcal{I}\mathcal{G}_{\Delta t}\mathcal{R}$ being a good approximation of $\mathcal{F}_{\Delta t}\mathcal{IR}$.
However, in heterogeneous wave media, $\mathcal{F}_{\Delta t} \mathcal{IR} \mathbf{u}$ will likely produce higher frequency modes, even if $\mathcal{R}\mathbf{u}$ have only low frequency modes. Consequently, we should not expect the approximation by $\mathcal{I}\mathcal{G}_{\Delta t}\mathcal{R}$ to be globally  sufficient. 


One idea is to replace $\mathcal{I}$ in $\mathcal{I}\mathcal{G}_{\Delta t}\mathcal{R}$ by something more elaborate. 
Our objective is to 
construct neural networks that approximate the operation 
\begin{equation}
    \mathcal{\theta}_{\Delta t}(\mathbf{w}, c)\approx \mathcal{F}_{\Delta t}\mathcal{I}\mathcal{G}_{\Delta t}^{-1}\mathbf{w}.
\end{equation}
The neural networks will be trained by data of the form $(c, \mathcal{G}_{\Delta t} \mathcal{R} \mathbf{u}, \mathcal{F}_{\Delta t} \mathbf{u})$, with $\mathbf{u}$ and $c$ sampled from certain distributions of interests.  
We then use the resulting enhanced solver
\begin{equation}\label{eq:theta-G-R}
    \theta_{\Delta t}\mathcal{G}_{\Delta t}\mathcal{R}\equiv \theta_{\Delta t}\left(\mathcal{G}_{\Delta t}\mathcal{R}\right)
\end{equation}
for wave propagation on the fine grid. Note that in most parts of this paper the solvers' dependence on the wave speed, $c$, is hidden for brevity of notation. 

We will rely on the parareal iteration to add back the missing high frequency components. The enhanced parareal scheme take the same form as above
\begin{equation}
\begin{aligned}
    \mathbf{u}_{n+1}^{k+1}:=  \theta_{\Delta t}\left(\mathcal{G}_{\Delta t}\mathcal{R}\right) \mathbf{u}_n^{k+1} + \left[\mathcal{F}_{\Delta t} \mathbf{u}_n^{k} - \theta_{\Delta t}\left(\mathcal{G}_{\Delta t}\mathcal{R}\mathbf{u}_n^{k}\right)\right],~~k=0,1,\dots, K,
\end{aligned}
\end{equation}
with
\begin{equation}\label{eq:u0_n-enhanced-stepping}
    \mathbf{u}_{n+1}^0:=\theta_{\Delta t}\left(\mathcal{G}_{\Delta t} \mathcal{R}\mathbf{u}_n^0\right).
\end{equation}

The wave solution operator, and most numerical solvers, are linear and reversible with respect to the initial data but  nonlinear to the wave speed. So is $\theta_{\Delta t}.$
However, the operator $\theta_{\Delta t}(\mathbf{u},c)$ is a so-called "near-identity" operator for sufficiently smooth wave speeds and wave fields.  Therefore, nearby a fixed reference speed function $c_0$, the linear operator $\theta_{\Delta t}(\cdot,c_0)$ has the potential to offer an adequate, though non-optimal, approximation of $\mathcal{F}_{\Delta t}\mathcal{I}(\mathcal{G}_{\Delta t})^{-1}$. 

It may be helpful to inspect the approximation of an individual Fourier mode of the wave field and consider errors in the phases and the amplitudes.
The additive correction in the parareal iterations is very effective in reducing the amplitude errors. However, the phase errors generally have a detrimental effect on the stability of parareal iterations. See e.g. \cite{ariel2016parareal}\cite{mikioInfluence18}\cite{nguyen2020stable}\cite{Ruprecht18}. 
It is therefore imperative that $\theta_{\Delta t}$ corrects the numerical dispersion errors in $\mathcal{G}_{\Delta t} \mathcal{R}\mathbf{u}$ compared to $\mathcal{F}_{\Delta t}\mathbf{u}.$ Furthermore, as discussed above, we also expect that $\theta_{\Delta t}$ to introduce some of the higher frequency modes that come out of propagating a smooth given wave field through the given non-constant media.  

Notice that for each instance of $c$, $\mathcal{F}_{\Delta t}\mathcal{I}\mathcal{G}_{\Delta t}^{-1}$ is a linear operator, provided that the prolongation $\mathcal{I}$ is linear.
With the nonlinear activations and biases in the network removed, the optimization model defines a regression model to construct data-driven (factorized) linear approximations of the operator associated with $\mathcal{F}_{\Delta t}\mathcal{I}\mathcal{G}_{\Delta t}^{-1}$. 
In which regime will this type of linear approximation work?
We shall present some experiments and comparisons regarding this observation. 

In Section~\ref{sec:J-net}, we present the proposed learning model and our approach in generating training data.
In Section~\ref{sec:effective_propagation}, we report results showing that a properly trained, simple, multi-level network can be an effective corrector of the coarse solver in a range of media.  
In Section~\ref{sec:parareal}, we apply the proposed networks to the parareal scheme in four wave speed models. 

\subsection{Other deep learning related work}

The authors in \cite{siahkoohi2019neural}\cite{rizzuti2019learned} train a convolutional neural network (CNN) for correction of low-fidelity solutions.  
As a supervised learning task, the data consists of low-fidelity and high-fidelity pairs, computed using the second-order and twentieth-order finite-difference, respectively. The data are time snapshots of the solutions for various media and sources. The experimental results in these papers indicate that the trained CNN can remove the numerical dispersion in complex velocity models. This result is encouraging since their setup is quite straightforward. 

In \cite{moseley2020solving}, a neural network model is proposed to learn wave propagation on a fixed medium. The Physics-Informed Neural Network (PINN) \cite{raissi2019physics} approach to setup the neural network's optimization model. This means the use of the given PDE on a chosen setup of points in the domain as the penalty term. The fully-connected neural network learns the linear relation between the initial wave source and the resulting wave field at later time. Each forward feed outputs the estimated wave field at the specified location. Applying the network for simulation on a different wave medium requires retraining the network. 
\LL{Ovadia \textit{et. al.}  \cite{ovadia2021beyond} proposed a PINN-based method 
for the homogeneous acoustic wave equation in one dimension. The resulting explicit nonlinear scheme seems stable for a step size that formally violates the Courant-Friedrichs-Levy condition. }

Another approach to using neural networks in a parareal framework is discussed in \cite{meng2019ppinn}. The authors propose parareal physics-informed neural network (PPINN). In particular, the sequence of fine solvers, mapping solutions from $t_n$ to $t_{n+1}$, $n=1,2,...$, are fully connected neural networks  (PINNs) trained with
loss functions that incorporate penalize the network outputs misfit to given differential equations. See  \cite{raissi2019physics}.  
The key idea of \cite{meng2019ppinn} is using the the parareal scheme to train the sequence of PINNs in parallel.
However, to the best of our current knowledge, the accuracy and stability of the PINN method for wave problems are not well understood. 

\LL{Recently,  Kochkov \textit{et. al.} \cite{kochkov2021machine} proposed a machine learning approach to speed up fluid dynamics simulations.
The approach enhances the accuracy of computations on the coarse grids. 
The ground truth data come from highly accurate computations performed on finer grids. 
Our philosophy is similar to theirs in that we propose to use machine learning to improve the classical numerical methods.}
Our approach, initiated in \cite{nguyen2020thesis}, also shares a certain similarity to \cite{siahkoohi2019neural}, but with additional key findings and further improvements. Most differently, we use the wave energy variables to define the optimization problem in the supervised learning setup; i.e. the spatial gradient $\nabla u$ and the time derivative of $u$ weighted by the wave speed. 
We think that this strategy results in a more convex optimization formulation and allows the algorithm to retain more physical properties of the wave system. This setup also leads to the input variables being dimensionless, which allows the algorithm to scale better. Such approach involving the energetic variables are used successfully in e.g. \cite{nguyen2020stable}\cite{rocha2018elastic}\cite{ rocha20193d}\cite{rocha2016acoustic}\cite{Tanushev-Engquist-Tsai}.
\LL{Finally, we advocate using parareal iterations to further improve the machine-learning based enhancement while providing a parallelism in time evolution.}

\subsection{The numerical solvers}\label{sec:numerical-solver}
In this section, we describe the coarse and fine solvers used in this paper.  \LL{Again, we point out that the choice of the fine and coarse solvers can be problem dependent. The main point of this paper is to introduce a deep learning approach to improve the fidelity of the coarse solver, and consequently enables stable parallel-in-time computation of the fine solver.}

The time stepping of both coarse and fine solvers are built as follows. Starting with  
\begin{equation}\label{eq:discretoize-Laplace-u}
    \partial_{tt} u(x,t) \approx
    c(x)^{2} Q_h u(x,t),
\end{equation}
where $Q_h u(x,t)$ denotes a numerical approximation of $\Delta u(x,t).$

The central scheme in time, which is the Leap Frog scheme for \eqref{eq:discretoize-Laplace-u}, can be written as the velocity Verlet time integrator for $\mathbf{u}:=(u,u_t)$: (with $v$ approximating $u_t$) 
\begin{equation}\label{eq:central_scheme_Verlet}
    \begin{aligned}
    u(x,t+k) =& u(x,t) + k v(x,t) + k\lambda~c(x)^{2}Q_h u(x,t)\\
    v(x,t+k) =& v(x,t) +\lambda\left( Q_h u(x,t)+Q_h u(x,t+k)\right)
    \end{aligned}
\end{equation}
where $\lambda=k/2h^2$. 
We will denote the mapping defined in \eqref{eq:central_scheme_Verlet} by
$\mathcal{S}^{Q_h}_{h,k}(\mathbf{u}),$ hiding the dependence on $Q_h$ for the brevity of notation.

\paragraph{The coarse solvers} 
$\mathcal{G}_{\Delta t^*}$, is defined by running the stepping scheme $\mathcal{S}^{Q_h}_{\Delta x,\Delta t}$ for $M$ time steps on the coarse grid $\Delta x\mathbb{Z}^2\times \Delta t\mathbb{Z}^+$; i.e.
\[
\mathcal{G}_{\Delta t^*}:=(\mathcal{S}^{Q_h}_{\Delta x,\Delta t})^M,~~~\Delta t^*=M\Delta t.
\]
$Q_h$ 
 is defined by the standard second order central differencing for the partial derivatives, applied on uniform Cartesian grids for the spatial domain.



\paragraph{The fine solver} $\mathcal{F}_{\Delta t}$ that we consider in this paper is formally denoted as 
\[\mathcal{F}_{\Delta t^*}:=(\mathcal{S}^{Q_h}_{\delta x, \delta t})^m,~\Delta t^*=m\delta t,\]
\LL{where $Q_h$ is either the second order central differencing or a spectral approximation of $\Delta u$, on a finer grid: $\delta x\mathbb{Z}^2\times \delta t \mathbb{Z}^+$.  
Some of the examples presented below involve discontinuous wave speeds, where the central scheme is used.
In Section~\ref{sec:unscaled-Marmousi-sim}, we present results involves the  pseudo-spectral method as the fine solver, since it is also used in the seismic imaging community.}

We assume that both the fine and the coarse solvers are stable on the respective grids.
Furthermore,  we assume that the fine solver is sufficiently accurate for wave speed $c\in(0,c_{max}]$, where $c_{max}$ is an upper bound for the wave speed in the class of wave media of interest.

If the two solvers are defined on different Cartesian grids ($\Delta x > \delta x$), we need to employ additional operators to map grid functions defined on the two grids. 
In this paper, $\mathcal{R}$ denotes the restriction operator that maps a finer grid function defined on $\Delta x\mathbb{Z}^2$ to one on $\delta x\mathbb{Z}^2$. 
$\mathcal{I}$ denotes the prolongation operator that maps a coarse grid function to fine grid. 

With a slight abuse of notations, we will also use $u$, $u_t$, and $\mathbf{u}=(u,u_t)$ to denote the corresponding discrete versions computed on the fine grid.

\paragraph{The wave energy semi-norm}
\LL{The wave equation is well-posed in the energy semi-norm}
\begin{equation}
    E[\mathbf{u}]:=\frac{1}{2}\int_{[-1,1]^2} |\nabla u|^2 + c^{-2}|u_t|^2 dx.
\end{equation}
\LL{This means, in particular, that wave propagation is stable in the wave energy semi-norm with respect to perturbation to the wave field. 
Furthermore, due to linearity, if $u$ and $v$ are solutions of the wave equation, $u-v$ also satisfies the wave equation. Wave fields with zero energy are constants and do not propagate. Therefore, for the purpose of wave propagation, it is natural to compare wave fields using the energy semi-norm. In this paper, we shall use the discretized version, which defines a discrete semi-norm on the grid:}
\begin{equation}\label{eq:discrete_energy_norm}
    E_h[\mathbf{u}] := \sum_{x_{i,j}\in h\mathbb{Z}\cap[-1,1]^2} \Big( \|\nabla_h u(x_{i,j})\|^2_2 + |c^{-1}(x_{i,j})\partial_t u(x_{i, j})|^2 \Big) h^2.
\end{equation}

In this paper, we will design algorithms that are defined as functions of the energy components $(\nabla u,u_t)$ of the wave field $\mathbf{u}=(u, u_t)$. \LL{In the following, we shall argue that the correction of numerical dispersion errors can be conveniently carried out by convolutions of the wave energy components of the given wave field. }

On the spatial grid $h\mathbb{Z}^2\cap[-1,1]^2$ with periodic boundary conditions,
define the mapping into energy components as
\begin{equation}
    \Lambda_h:\mathbf{u}\mapsto(\nabla_h u, c^{-1} u_t),
\end{equation}
where $\nabla_h u$ is computed by numerically. $\Lambda^\dag_h$  denotes the pseudo-inverse of $\Lambda_h$ and is defined in Appendix~\ref{Appendix:Lambda_Lambda_inv}. 

In \cite{rocha2018elastic, rocha20193d, rocha2016acoustic}, the authors use the energy norm in seismic imaging methods and their results show improved reconstructions over the regular $\ell_2$ loss. The Gaussian beam method \cite{liu2013error}\cite{Tanushev-Engquist-Tsai} also utilizes the energy norm to decompose wavefields. Following the success of these work, our deep learning setup will base on the energy components induced by the energy norm. \LL{See Section~\ref{sec:netexample}.}

\paragraph{The enhanced solvers}
Our goal is to construct neural network functions $\mathcal{H}_{w}$, parameterized by $w$, such that  
$$
\mathcal{H}_{w}(\Lambda_{\Delta x}~{\mathcal{G}}_{\Delta t^*} {\mathcal{R}}\mathbf{u}, c) \approx  \Lambda_{\delta x} \mathcal{F}_{\Delta t^*} \mathbf{u}
$$
For a class of piecewise smooth wave speed functions. In other words, $\mathcal{H}_{w}$ corrects the energy components computed by the coarse solver $\mathcal{G}_{\Delta t^*}\mathcal{R}$,  according to those computed by $\mathcal{F}_{\Delta t^*}$. 

With $\mathcal{H}_w$, we define the enhanced solver as
\begin{equation}
    \theta_{\Delta t^*}(\mathcal{G}_{\Delta t^*}\mathcal{R}\mathbf{u}, c):= \Lambda^\dag_{\delta x}\mathcal{H}_w(\Lambda_{\Delta x} \mathbf{u}, c).
\end{equation}
\LL{The enhanced solver, with step size $\Delta t^*$, operates in a regime beyond the limitation from the Courant-Friedrichs-Levy (CFL) condition for linear schemes, propagating information across a distance of several $\delta x$ in each time step. }


\begin{table}
    \centering
    \caption{A list of key notations.}
    \begin{tabular}{l|l}
       \hline
       $\Delta t^*$ & The time step size of the enhanced coarse solver\\
       $\textbf{u}$ & Wave field vector: $(u, \partial_t u)$\\
       $\theta_{\Delta t^*}$ & Correction operator\\    
       $\mathcal{G}_{\Delta t^*}$ & Coarse solver defined on a coarse grid\\
       $\mathcal{F}_{\Delta t^*}$ & Fine solver\\
       $\theta_{\Delta t^*}\mathcal{G}_{\Delta t^*}\mathcal{R}$ & The enhanced solver operating on the fine grid \\
       $\mathcal{H}_w$  &  A neural network used for defining $\theta_{\Delta t^*}$ \\
       $E_h[\textbf{u}]$ & The discrete wave energy norm of $\mathbf{u}$.\\
       $\mathtt{D_t}(\Delta t^*)$ & An elementary training data set\\
       $\mathtt{D^p_t}(\Delta t^*)$ & The training data set derived from $\mathtt{D_t}$\\
       $\mathsf{JNet(L,s)}$ & \LL{A network with $L$-level encoder and $(L+s)$-level decoder} \\
       \hline
    \end{tabular}
    \label{tab:notations}
\end{table}

\subsubsection*{Correction of dispersion errors}\label{sec:correction_op}

Here we present a simple and classical dispersion property of the central scheme. The point is to identify a regime in which $\theta_{\Delta t^*}(\mathbf{u}, c)$  is an identity plus a minor perturbation. In such a regime, the enhanced solver should approximate the perturbation.

We analyze a simplified setup involving plane waves: $\exp{i(\omega t+kx)}$ propagating in a medium with constant speed $c$ in one dimension. 
We assume that the fine solution operator $\mathcal{F}_{\Delta t}$ to be exact, mapping the given wave field at time $t$ to $t+\Delta t^*$:
$$
    \mathbf{u}(x,t+\Delta t^*) = \mathcal{F}_{\Delta t^*} \mathbf{u}(x,t),
$$
which has the linear dispersion relation $\omega(c,k) = \pm c k.$ Here we exclusively use $k$ to denote the wave number (not the parareal iteration). 
On the other hand, we approximate the fully-discrete coarse solver
by a semi-discrete scheme. More precisely, in this scheme denoted as $\mathcal{G}_{\Delta t^*}^{sd}$, central differencing  with step size $\Delta x$ is used to approximate the spatial derivatives, and the resulting system is integrated in time exactly. 
We denote the coarse solution

$$
    \Tilde{\mathbf{u}}(x,t+\Delta t^*) = \mathcal{G}_{\Delta t^*}^{sd} \mathbf{u} (x,t). 
$$
The semi-discrete scheme $\mathcal{G}^{sd}_{\Delta t^*}$ gives us the numerical dispersion relation 
\begin{equation}\label{eq:num_dispersion_expression}
{\omega}_{\Delta x}(c,k) = \pm c k \Big(1  -  (\dfrac{k^2\Delta x^2}{24}- \dfrac{k^4 \Delta x^4}{1920} + \mathcal{O}(k\Delta x)^6))\Big).
\end{equation}
Define the difference between the {positive-signed} $\omega$ and $\omega_{\Delta x}$ as the numerical dispersion error 
\begin{equation}\label{eq:num_dispersion_err}
     \varepsilon_{c,\Delta x}(k)= ck \Big(\dfrac{k^2\Delta x^2}{24}- \dfrac{k^4 \Delta x^4}{1920} + \mathcal{O}(k\Delta x)^6 \Big).
\end{equation}
We see the well-known property that higher wave numbers lead to larger dispersion errors, and asymptotically, the need to decrease $\Delta x$ at a faster rate than $1/k$. 

We can thus approximate the energy components of the fine solution, $(\nabla u, c^{-1} u_t)$ by correcting the numerical dispersion errors in the coarse solutions $(\nabla \tilde u, c^{-1} \tilde u_t)$. The correction can be written as a convolution: {
\begin{equation} \label{eq:omegaconv}
    \left[\begin{array}{cc}
        \nabla u \\
        c^{-1} u_t
    \end{array}\right]
    \simeq
    \left[\begin{array}{cc}
        \mathscr{F}[\cos({\varepsilon_{c,\Delta x} ~\Delta t^*})] & \pm \mathscr{F}[\sin({\varepsilon_{c,\Delta x} ~\Delta t^*})] \\
        \mp \mathscr{F}[\sin({\varepsilon_{c,\Delta x} ~\Delta t^*})] & \mathscr{F}[\cos({\varepsilon_{c,\Delta x} ~\Delta t^*})]
    \end{array}\right] \ast
    \left[\begin{array}{c}
        \nabla \Tilde{u} \\
        c^{-1}\Tilde{u}_t
    \end{array}\right],
    \end{equation}}
where $\mathscr{F}[\cdot]$ denotes the Fourier transform from the frequency $k$ to $x$. 
For details, readers refer to Appendix~\ref{Appendix:analytic_correction_op}.

We notice that in Equation~\eqref{eq:num_dispersion_expression}, 
for small $k\Delta x, \Delta t^*$, the leading order term in the correction is an identity map.  This fact could serve as a theoretical motivation for adding a skip connection between every corresponding pair of layers in the multi-level network described below.
The numerical dispersion relation \eqref{eq:num_dispersion_expression} and \eqref{eq:omegaconv} also predict that the main role of the neural network would be to overcome the challenges resulting from higher frequency wave components, to say the least.

{In Section~\ref{sec:J3-J6-dispersion-comp}, we presents an example illustrating the numerical dispersion errors in a piecewise constant medium. There, we will compare the correction of such errors by the proposed networks.}

\section{The supervised deep learning setup}\label{sec:J-net}

In this section, we describe the network architecture, the optimization model, and the generation of training examples.

The network's input variables are the energy components of a wave field on the coarse grid and the wave speed function --- these are functions on teh coarse grid.
The outputs are the energy components of the processed wave field, but
on the fine grid.
The network does not explicitly involve any approximation of the derivatives of the wave field. 

\LL{Since waves have finite speed of propagation, it suffices to learning only locally in time and space how a family of wave fields propagate through a family of representative wave media in the $\Delta t^*$ time interval. Since the networks are based on convolutions using kernels with compact width, they may be applied to larger input tensors. }

\subsection{The \texttt{JNet} architecture}

Our chosen network architecture is a U-net architecture \cite{ronneberger2015unet} with skip connections, though in our setup the downward-upward path in the network resembles more like the letter J.

\paragraph{Multiple levels with convolution layers.}
The network contains a cascade of convolution layers encoding the input tensor into a so-called feature space. Then a cascade of convolution layers decodes the feature space to produce the final output image. The encoding pathway utilizes average-pooling layers to downsample the intermediate output, while the decoding pathway uses the bilinear interpolation. 

\begin{figure}
    \centering
    \includegraphics[width=1.0\linewidth]{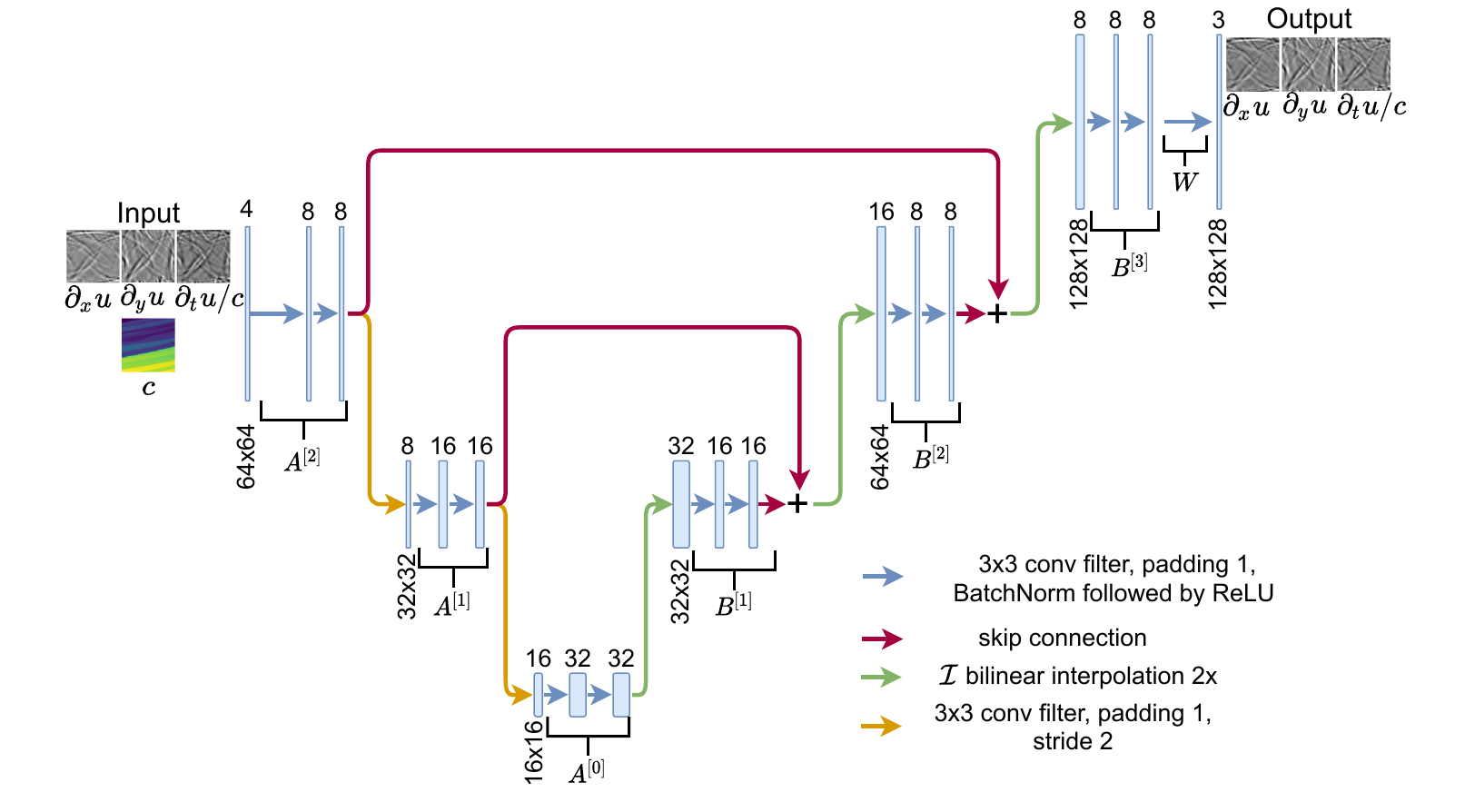}
    \caption{Diagram of a $\mathsf{JNet(3,1)}$. If the network uses ReLU activations, then there are bias parameters in each conv filter, followed by a ReLU activation and a batch norm layer.}
    \label{fig:Jnet}
\end{figure}

\paragraph{Notations. }\LL{We introduce the notation $\mathsf{JNet(L,s)}$ for $L$-level encoder, and $(L+s)$-level decoder. A $\mathsf{JNet(L,0)}$ is an $L$-level U-net.}  
Unless mentioned otherwise, the networks are activated by the Rectified Linear Units (ReLUs). We call a  $\mathsf{JNet(L,s)}$ as a linear, when the ReLU activation is replaced by the identity function and no biases removed. Throughout the paper, a $\mathsf{JNet(L,1)}$ network is also referred to as an $L$-level (linear or ReLU) network. 
Figure \ref{fig:Jnet} illustrates a $\mathsf{JNet(3,1)}.$

\begin{figure}
    \centering
    \includegraphics[width=1.0\linewidth]{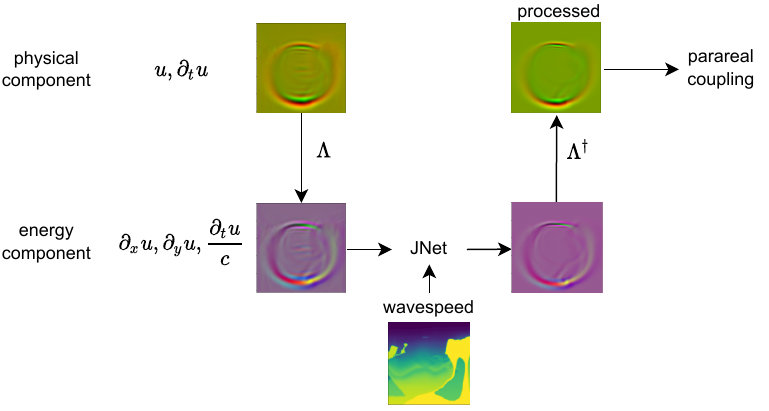}
    \caption{Diagram of how a $\mathsf{JNet}$ enhanced the wave propagation.}
    \label{fig:DNN_parareal}
\end{figure}

\paragraph{The input tensor.} The input of the network is a tensor that has $4$ channels, corresponding to the energy components of the coarse solution and the wave speed functions defined on the coarse grid.

In particular, in the remaining of the paper, except for Section~\ref{sec:unscaled-Marmousi-sim}, we shall report computational results obtained from trained $\mathsf{JNet(3,1)}$ and $\mathsf{JNet(6,1)}.$ 
In this paper, we will compare these networks after training with different databases. 

Finally, as a reference, the forward feed of a $\mathsf{JNet(3,1)}$, implemented in PyTorch with a $128\times 128\times 4$ input tensor, takes roughly $0.028$ seconds on the CPU of a moderately equipped laptop. 
On the same computer, the fine solver, implemented in NumPy, takes roughly $0.058$ seconds to evolve a given initial wave field of size $128\times 128\times 2$ to time $t=0.1.$

\subsection{Data sets and optimization model}\label{sec:netexample}

To learn the (parametric) coarse-to-fine solution map
$
    \mathcal{H}_{w}:X \mapsto Y, 
$ 
the training data should contain sufficiently representative examples, sampling the wave speed $c$ of interest and the wave fields of relevance --- those propagated by the coarse and fine solvers, starting from a class of initial conditions of interest.  
In our setup, 
\begin{eqnarray*}
   X &:= &\{\Lambda_{\Delta x} \mathcal{G}_{\Delta t^*} \mathcal{R}\mathbf{u}, c \}, \\
   Y &:= &\{\Lambda_{\delta x} \mathcal{F}_{\Delta t^*} \mathbf{u}\}. 
\end{eqnarray*}
The wave speed $c$ and the wave field $\mathbf{u}$ are sampled from the distributions described below. 

Both the coarse and the fine solvers use the standard second order scheme defined in Section~\ref{sec:numerical-solver}.
The discretization parameters of the solvers are tabulated in Table \ref{tab:disparam}. 

The data sets discussed below will be made available for download.

\begin{table}
\caption{Discretization parameters for the fine and coarse solver.}
\label{tab:disparam}
\medskip
\centering
\begin{tabular}{cccccc} 
	\hline
     $\Delta t^*$ & $\delta x$ & $\delta t$ & $\Delta x$ & $\Delta t$ & order of interp  \\
     \hline
     $0.1, 0.2, 0.25$ & $2/128$ & $1/1280$ & $2/64$ & $1/160$ & $4$-th order \\
	\hline
\end{tabular}
\end{table}

\subsubsection*{The wave speed and initial wave field samples $\mathtt{D_0}$ and $\mathtt{D_0^p}(\Delta t^*)$}

\begin{figure}
    \centering
    \includegraphics[width=4in]{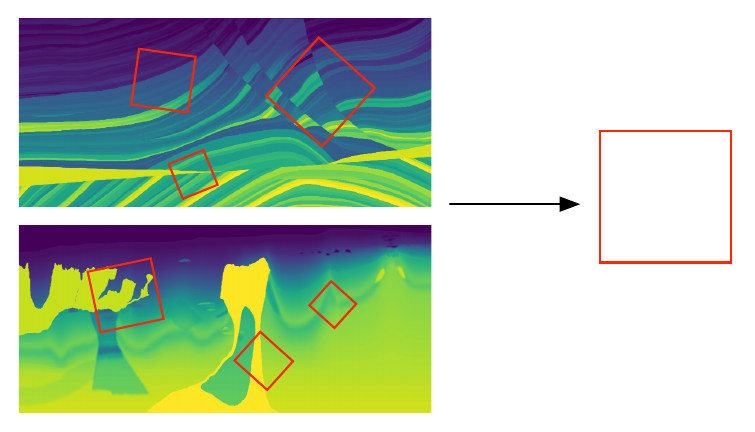}
    \caption{Square subregions of different sizes in the Marmousi profile are rotated and rescaled to a reference domain with the sample wave speed $c_s$. $\mathcal{F}_{\Delta t^*}$ and $\mathcal{G}_{\Delta t^*}$ are applied to different initial pulses on it to generate training data.}
    \label{fig:refence_wave_domain}
\end{figure}

The data set for the medium wave speed comprises of randomly selected subregions (of different areas and rotations) of the full Marmousi model \cite{brougois1990marmousi} and the BP model \cite{billette2004bpmodel}. 
The wave speed in each of the subregions will be mapped (subsampled and rescaled\footnote{The wave speed in each medium example is rescaled by an integer, allowing the coarse solver to be stable with the predetermined $\Delta x$ and $\Delta t$.}) to a $128\times 128$ array, which is regarded as a grid function defined on $h\mathbb{Z}^2\cup [-1,1)^2$, $h=2/128$. See Figure~\ref{fig:refence_wave_domain}.

Additionally, we remark that Wu et al. \cite{fomel2020buildingseismicvelo} developed a comprehensive methodology to synthesize the mediums with realistic earth layers and faults. Their approach may be use to further enrich the data set.

\paragraph{$\mathtt{D_0}$:}
For each of the sampled subregions, 
an initial wave field, $\mathbf{u}_0=(u_0, \partial_t u_0))$ is sampled from the Gaussian pulses of the form:  
$$ u_0(x,y) = e^{-(x^2+y^2)/\sigma^2}, \partial_t u_0(x,y) = 0, ~~~ x,y\in \delta x\mathbb{Z}^2\cap [-1,1)^2,~~~1/\sigma^2 \sim \mathcal{N}(250,10).$$
The pair $(c, \mathbf{u}_0)$ thus sampled is collected in $\mathtt{D_0}.$ 

Next, we generate a set of training examples that have more complex initial wave fields.

\paragraph{$\mathtt{D_0^p}(\Delta t^*)$:}
For each $(c,\mathbf{u}_0)$ in $\mathtt{D_0}$,
we collect $\mathbf{u}^k_n$ that are computed by the Procrustes parareal scheme, as defined in \eqref{eq:Procrustes-parareal}-\eqref{eq:Procrustes-optimization}, for $n=0,1,\dots, N,$ and $k=0,1,2,3,4$. We shall refer to this collection as $\mathtt{D_0^p}(\Delta t^*)$.

\subsubsection*{The example wave fields $\mathtt{D_t}(\Delta t^*)$ and $\mathtt{D_t^p}(\Delta t^*)$}

\begin{figure}
    \centering
    \includegraphics[width=5in]{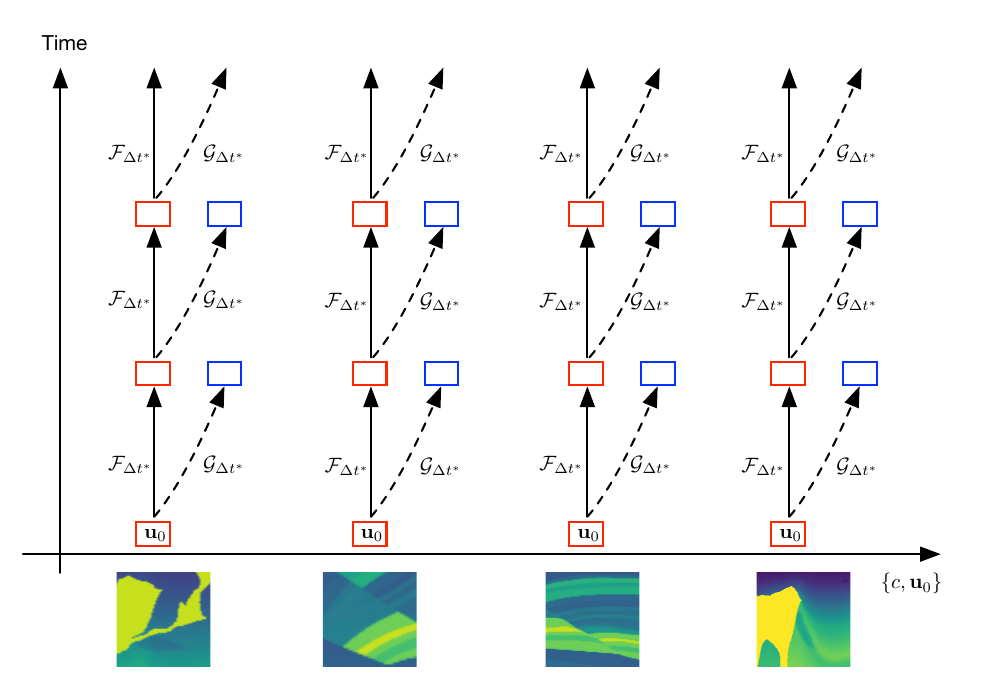}
    \caption{The sampling of the solution manifold:
On each random sample of medium wave speed,
the coarse and fine solvers propagate initial wave fields consisting of random Gaussian pulses. Of each data point, $x$ corresponds to a blue box above, and $y$ a corresponding red box.}
    \label{fig:sampling_wave_causality}
\end{figure}

We recognize the need to sample the strong causality in the wave dynamics. Ergo we sample the trajectories
passing through the wave speed-initial wave field pair in $\mathtt{D_0}$ and $\mathtt{D_0^p}.$ Figure~\ref{fig:sampling_wave_causality} demonstrate the scheme which we use to generate the training two training sets. 
The trajectories are discretized by a finite number of points computed from the fine solver (one could use the coarse solver for this as well). The points from the trajectories are paired with the associating wave speed function and  collected into the training data sets. 

\begin{description}
\item[$\mathtt{D_t}(\Delta t^*)$:]
We first collect the sequence of wave fields propagated by the fine solver: 
\[
\mathbf{u}_{n+1}:=\mathcal{F}_{\Delta t^*}\mathbf{u}_n,\,\,\, n=0,1,\dots,N-1,~~~\forall(c,\mathbf{u}_0)\in \mathtt{D_0}
\]
Then  
$\mathtt{D_t}(\Delta t^*)$ consists of the pairs
\begin{equation}\label{eq:from_D0_to_Dt}
\left(c, (\Lambda_{\Delta x}\mathcal{G}_{\Delta t^*} \mathcal{R}\mathbf{u}_n, \Lambda_{\delta x} \mathbf{u}_{n+1})\right),~~~ n=0,1,\dots,N-1.
\end{equation}

\item[$\mathtt{D_{t}^p}(\Delta t^*)$:]
This data set consists of the pairs
\begin{equation}
\left(c,(\Lambda_{\Delta x}\mathcal{G}_{\Delta t^*} \mathcal{R}\mathbf{u}, \Lambda_{\delta x}\mathcal{F}_{\Delta t^*} \mathbf{u})\right),~~~(c,\mathbf{u})\in \mathtt{D_0^p}.
\end{equation}

\end{description}

\begin{figure}
    \centering
    \includegraphics[width=1\linewidth]{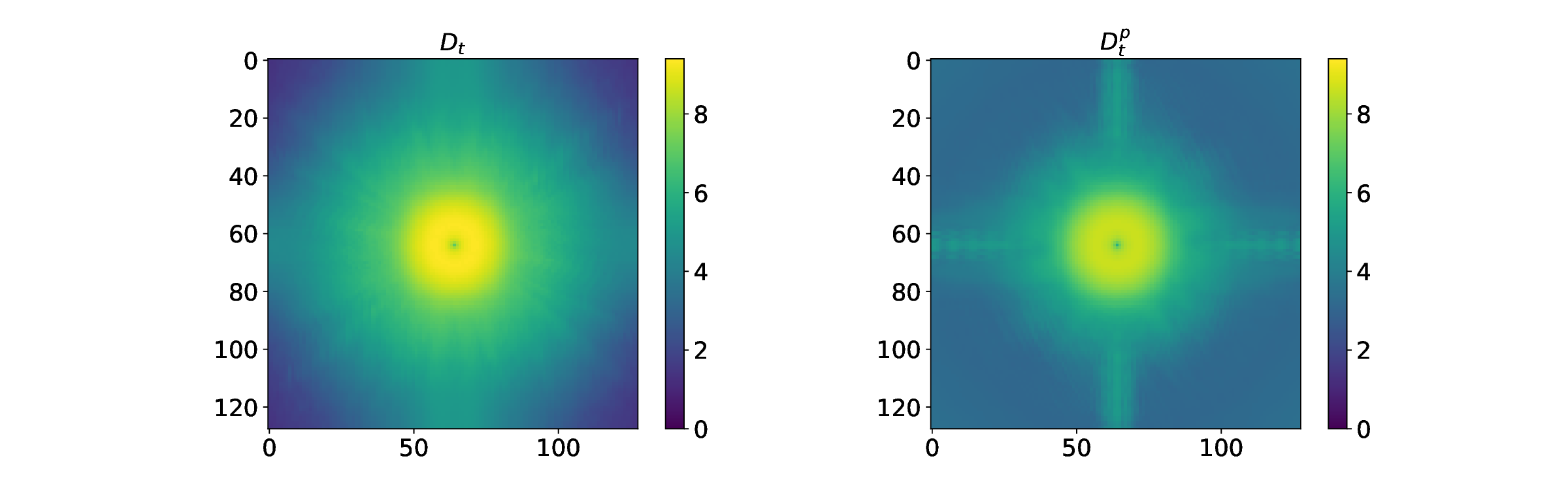}
    \caption{Sum of the squared Fourier modes magnitudes  of the $u_t$ component of $\mathcal{F}_{\Delta t^*}\mathbf{u}$ in the training wave field. The results are plotted in logarithmic scale. 
    }
    \label{fig:trainwavestats}
\end{figure}

We generate $10,000$ examples from $\mathtt{D_0}$ and $\mathtt{D_0^p}(\Delta t^*)$, and name the resulting databases 
$\mathtt{D_t}(\Delta t^*)$ and $\mathtt{D_t^p(\Delta t^*)}$ respectively. When the context is clear, we will use $\mathtt{D_t}$ and $\mathtt{D_t^p}$ for simplicity. 
{Figure~\ref{fig:trainwavestats} shows the Fourier mode statistics of the wave fields in $\mathtt{D_t}$ and $\mathtt{D_t^p}$. We observe that the wave fields in $\mathtt{D_0^p}$ have more high frequency modes.}

\paragraph{\LL{The optimization model} }


Formally, the loss function is the mean squared errors:
\begin{equation} \label{eq:lossdef}
    \min_{\mathcal{H}_w\in\mathsf{JNet}(L,s)} \mathcal{L} (w; D) = \min_{\mathcal{H}_w\in\mathsf{JNet}(L,s)} \dfrac{1}{|D|} \sum_{(c,(X,Y))\in D} \|\mathcal{H}_{w}(X) - {Y}\|^2_2,
\end{equation}
where $D$ is the training data set and $|D|$ denotes the number of training examples. $D$ is either $\mathtt{D_t}(\Delta t^*)$ or $\mathtt{D_t^p}(\Delta t^*)$ defined above.
{However,  the summands correspond to the approximation errors (the difference in the enhanced wave field and the "labeled" wave field) in the discrete wave energy semi-norm \eqref{eq:discrete_energy_norm}.} 
\LL{More precisely, writing $$ \mathcal{H}_w(X) -Y = \{(q_n(w),  p_n(w)):n=1,2,\cdots, N,\}$$
we have 
\begin{equation}
    \|\mathcal{H}_{w}(X) - {Y}\|^2_2 \equiv E_h[(q_n(w), p_n(w)],
\end{equation}
where $E_h$ is defined in \eqref{eq:discrete_energy_norm}, $q_n(w)$ corresponds to the resulting wave field difference, and $p_n(w)$ the time derivative of the wave field difference. 
}

\LL{Now, considering  $\mathtt{D_t}$ and the fact that it contains the propagation of the wave fields in $\mathtt{D_0}$ (see Equation \eqref{eq:from_D0_to_Dt}), we rearrange the summation in \eqref{eq:lossdef} into
\begin{equation}
    \mathcal{L}(w;D) = \frac{1}{Z_D} \sum_{(c,X,Y)\in \mathtt{D_0}}\left[\sum_{n=1}^{N}  \left(||\nabla q_n(w)||^2_2-\frac{1}{c^2}|| p_n(w)||^2_2 \right) + \frac{2}{c^2}||p_n(w)||^2_2\right],
\end{equation}
where $Z_D$ is the normalizing constant.}

\LL{Recall the Least Action Principle associated with the wave propagation:
$$ \min_u \int_0^T\int |\nabla u|^2 -\frac{1}{c^2} |u_t|^2~~dx dt. $$
In light of this, the proposed optimization model \eqref{eq:lossdef} can be interpreted as minimizing the sum of the discretized least action principle and the $L^2$ norm of the momentum of the mismatch (wave field).
The same arguments apply to the case of $\mathtt{D_t^p}(\Delta t^*)$. However, the precise consequence of this interpretation will be the subject of a future paper.} 

\LL{Similar strategies employing the wave energy to compare wave fields for wave propagation purposes have been used successfully, for example in \cite{ariel2011gaussian},\cite{nguyen2020stable},\cite{rocha2016acoustic},\cite{Tanushev-Engquist-Tsai}, and \cite{tanushev2011gaussian}.
In \cite{Tanushev-Engquist-Tsai}, it is shown numerically that in the wave energy semi-norm is convex with respect to operations such as translation and rotation of the wave fields which are well approximated by superposition of Gaussian Beams.
}

\begin{figure}
    \centering

    \includegraphics[width=0.47\linewidth]{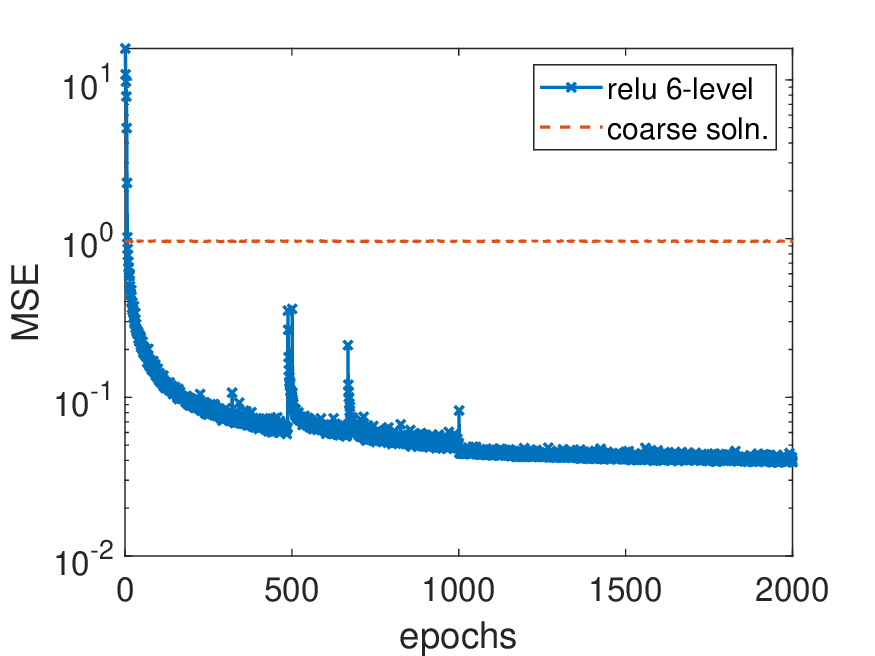}
    \includegraphics[width=0.47\linewidth]{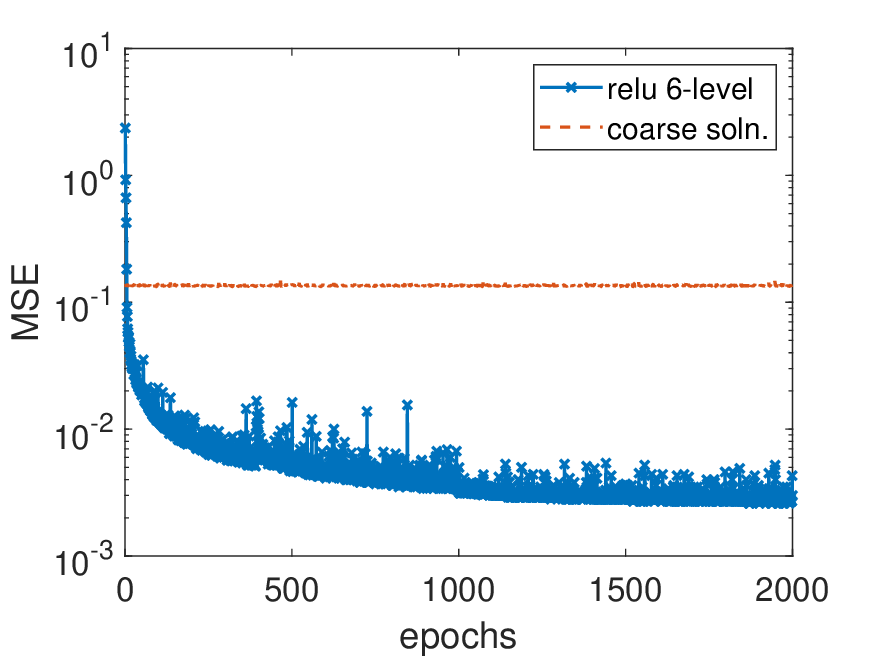}

    \caption{The training progress of two different data sets. Left: data set $\mathtt{D_t}$. Right: data set $\mathtt{D_t^p}$ containing snapshots of Procrustes parareal in crops of Marmousi and BP medium. Each data set has $10^4$ samples. In both training sessions, the batch size is 128, and the learning rate is $10^{-3}$.}
    \label{fig:trainingprogess}
\end{figure}

\paragraph{Training}
We use a simple stochastic gradient descent algorithm to minimize \eqref{eq:lossdef}.
The neural networks are trained by the mini-batch Adams, starting with the  initial conditions prescribed by PyTorch's built-in routines. 
The training data is randomly divided into batches of 128 examples. 
We report the MSE in a typical training progress for each data set in Figure~\ref{fig:trainingprogess}. 
Notice that the training errors on the right subplot appear approximately one order of magnitude smaller than those on the left. This is a subject of future investigation. It is possibly due to the Procrustes solutions's being closer to each other in some suitable sense. 

In Table~\ref{tab:train-test-errors}, we further report a comparison of the final training errors and testing errors between  the 3-level and the 6-level \texttt{JNet}s, trained  using $\mathtt{D_t^p}(\Delta t^*)$ with $\Delta t^*=0.2$ and $0.25$.

\begin{table}
\caption{Train and test errors of $\mathsf{JNet(3,1)}$ and $\mathsf{JNet(6,1)}$. The error is evaluated on the batch size of $128$. We trained these networks on $\mathtt{D_t^p}$ for $2000$ epochs; the learning rate in the first 1000 epochs is set to $10^{-3}$, then reduced to $5\cdot 10^{-4}$. }
    \label{tab:train-test-errors}
    \centering
    \begin{tabular}{lcccc} 
        \hline
        \multicolumn{1}{c}{} & \multicolumn{2}{c}{$\Delta t^* = 0.2$} & \multicolumn{2}{c}{$\Delta t^* = 0.25$} \\
        \hline
         {} & 3-level  & 6-level & 3-level & 6-level \\
         \hline
        Final training error $\times 10^{-4}$ & $5.2$ & $3.8$ & $7.1$ & $6.9$ \\
        \hline
        Test error $\times 10^{-4}$ & $6.8$ & $7.2$ & $9.8$ & $9.1$ \\
        \hline
    \end{tabular}
\end{table}

\section{Propagation}\label{sec:effective_propagation}

In this section, we study the enhanced propagator
$\theta_{\Delta t^*}$ in a few setups. 
We shall report the computed results  using the relative energy error
$$
    \dfrac{E_{\delta x}[\mathbf{u} -\mathbf{v}]}{E_{\delta x}[\mathbf{v}]},
$$ 
where $\mathbf{u}$ is the computed wave field and $\mathbf{v}$ is the reference solution.

\subsection{Comparison of propagation errors on constant media}

We test the enhanced solvers on constant media of different wave speeds, using initial wave fields consisting of pulses of various widths. 
The networks are trained on data sets of different cardinalities, sampled from the same distribution. 

The test initial wave field is defined as follows:
$$
    \mathbf{w}_{test}^\sigma = \mathcal{F}_{\Delta t^*} \circ \mathcal{F}_{\Delta t^*} \mathbf{w}_0^\sigma,~~~\mathbf{w}_0^\sigma = ( e^{-(x^2+y^2)/\sigma^2},0),
$$
where $\sigma$ is sampled uniformly $\sigma^{-1} \in [5,20]$. Larger valuers of $\sigma^{-1}$ correspond to wave fields whose Fourier modes for higher wave numbers have larger magnitudes (``sharper" wave fields). 

Coarse and fine solutions are propagated
in a constant medium $c\in [0.1,3.0]$. 
The errors in the energy is defined as
\begin{equation}
    e(\sigma, c) := \dfrac{E_{\delta x}[\mathcal{H}_{w}  ({\mathcal{G}}_{\Delta t^*}\mathcal{R} \mathbf{w}_{test}^\sigma,c) - \Lambda_{\delta x} \mathcal{F}_{\Delta t^*} \mathbf{w}_{test}^\sigma ]}{E_{\delta x}[\Lambda \mathcal{F}_{\Delta t^*} \mathbf{w}_{test}^\sigma]} \label{eq:generror}.
\end{equation}
We remark that the dependence of $c$ for $\mathcal{G}_{\Delta t^*}$ and $\mathcal{F}_{\Delta t^*}$ are hidden from the notation. 

\begin{figure}
    \centering
    \includegraphics[width=1.\linewidth]{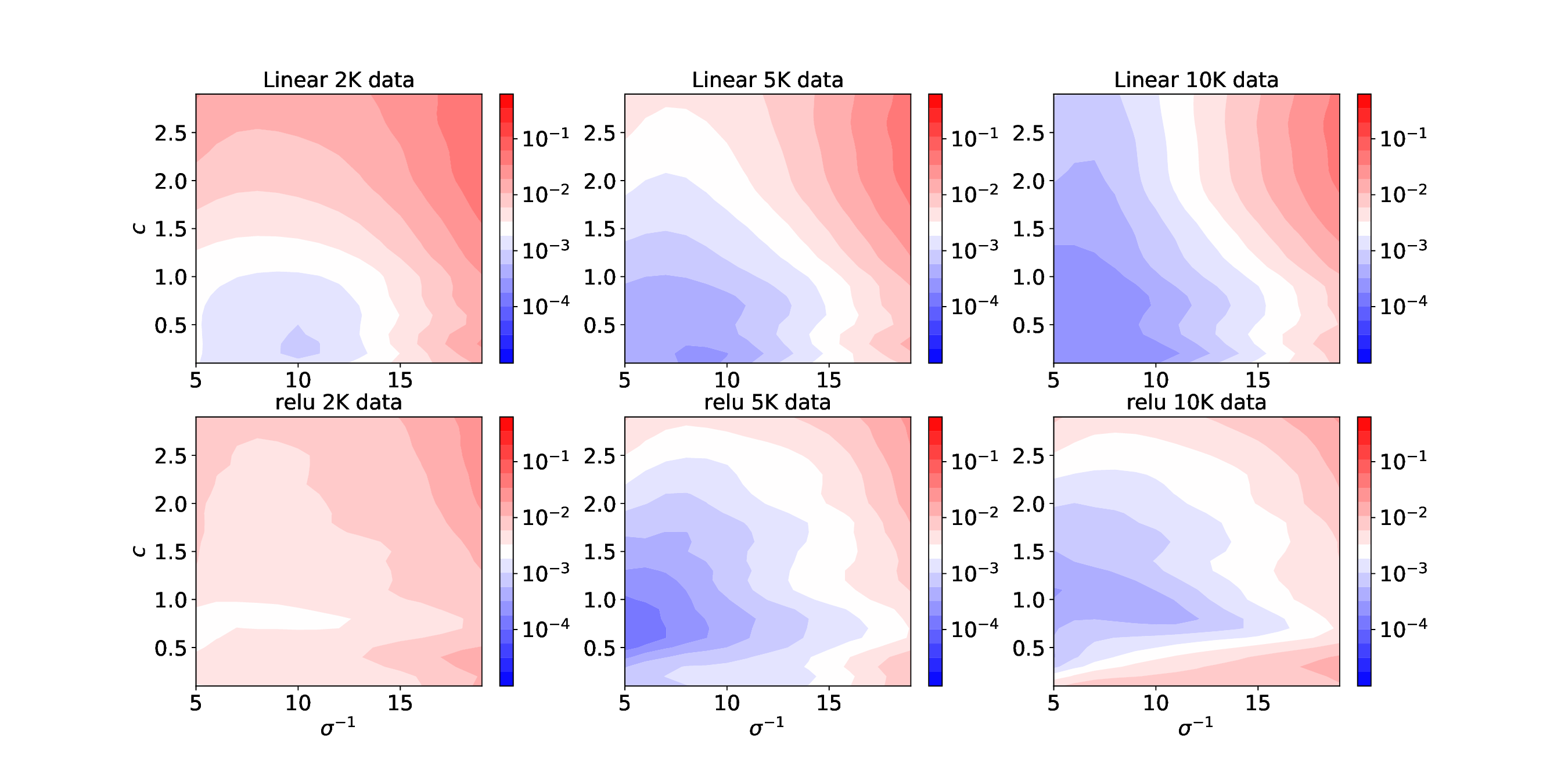}
    \includegraphics[width=0.6\linewidth]{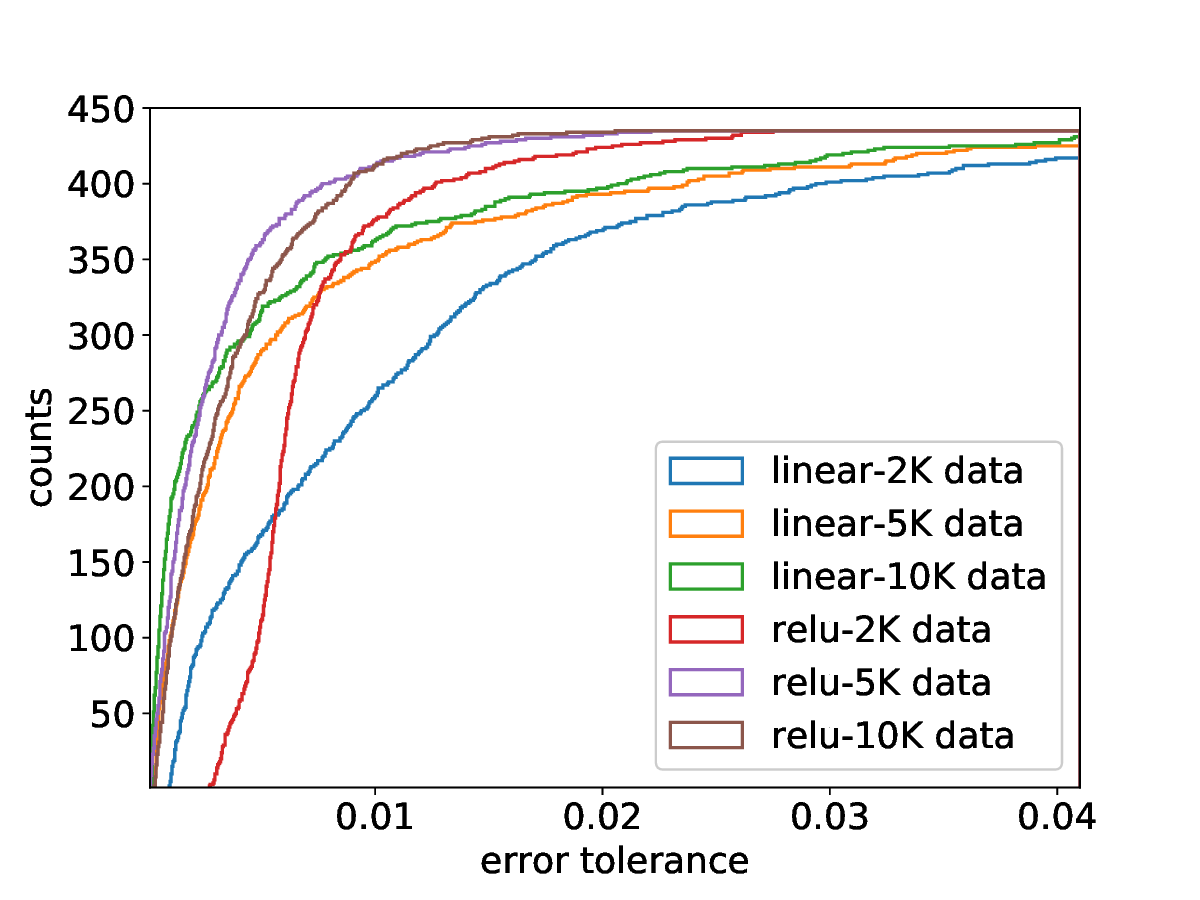}
    \caption{A comparison of the $e(\sigma, c)$ of linear and ReLU  $\mathsf{JNet(6,1)}$'s,  trained with different number of examples drawn from the same distribution, $\mathtt{D_t^p}(0.1)$. }
    \label{fig:proxy_error_comps}
\end{figure}

\begin{figure}
    \centering
    \includegraphics[width=1.\linewidth]{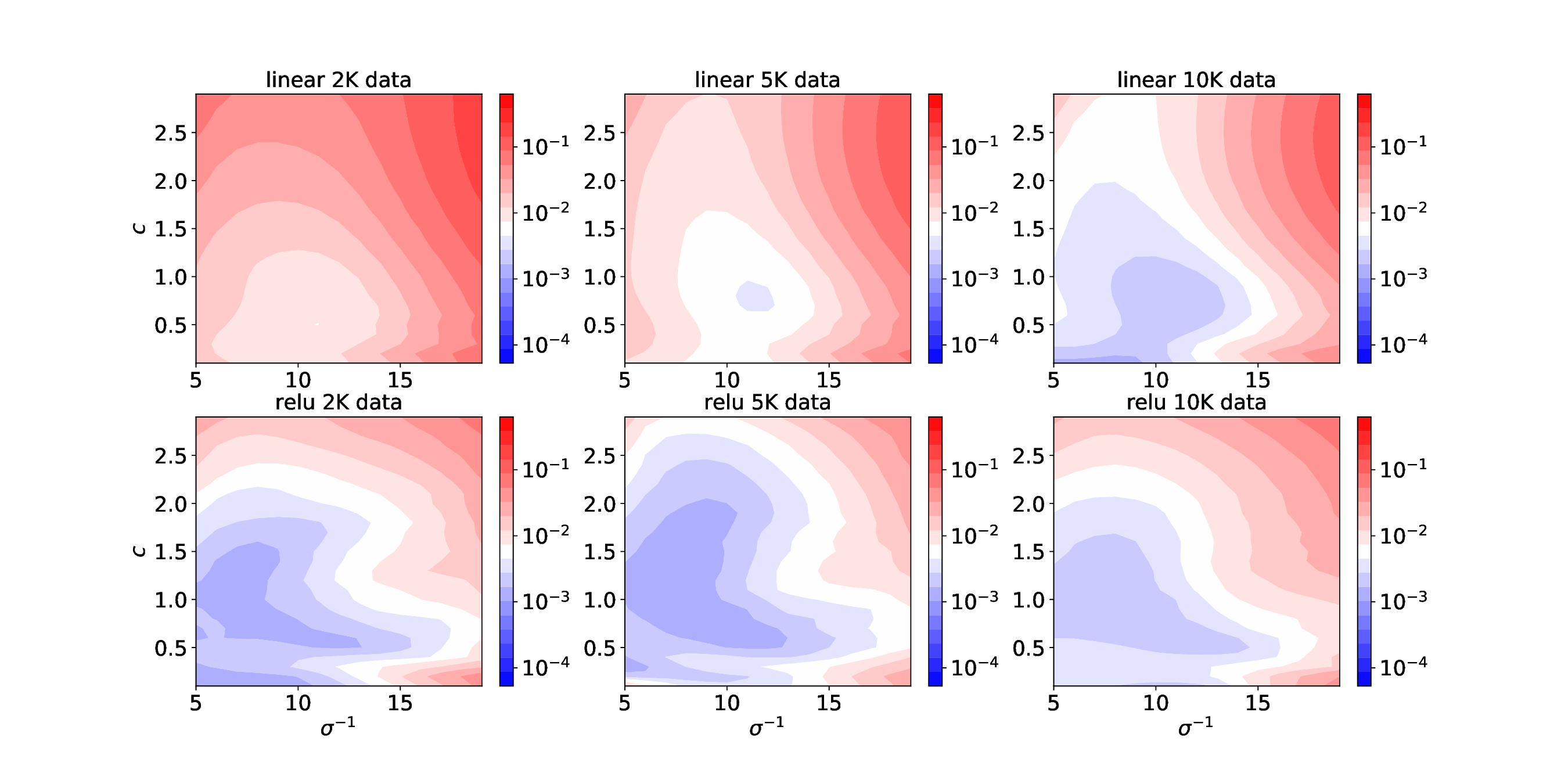}
    \includegraphics[width=0.6\linewidth]{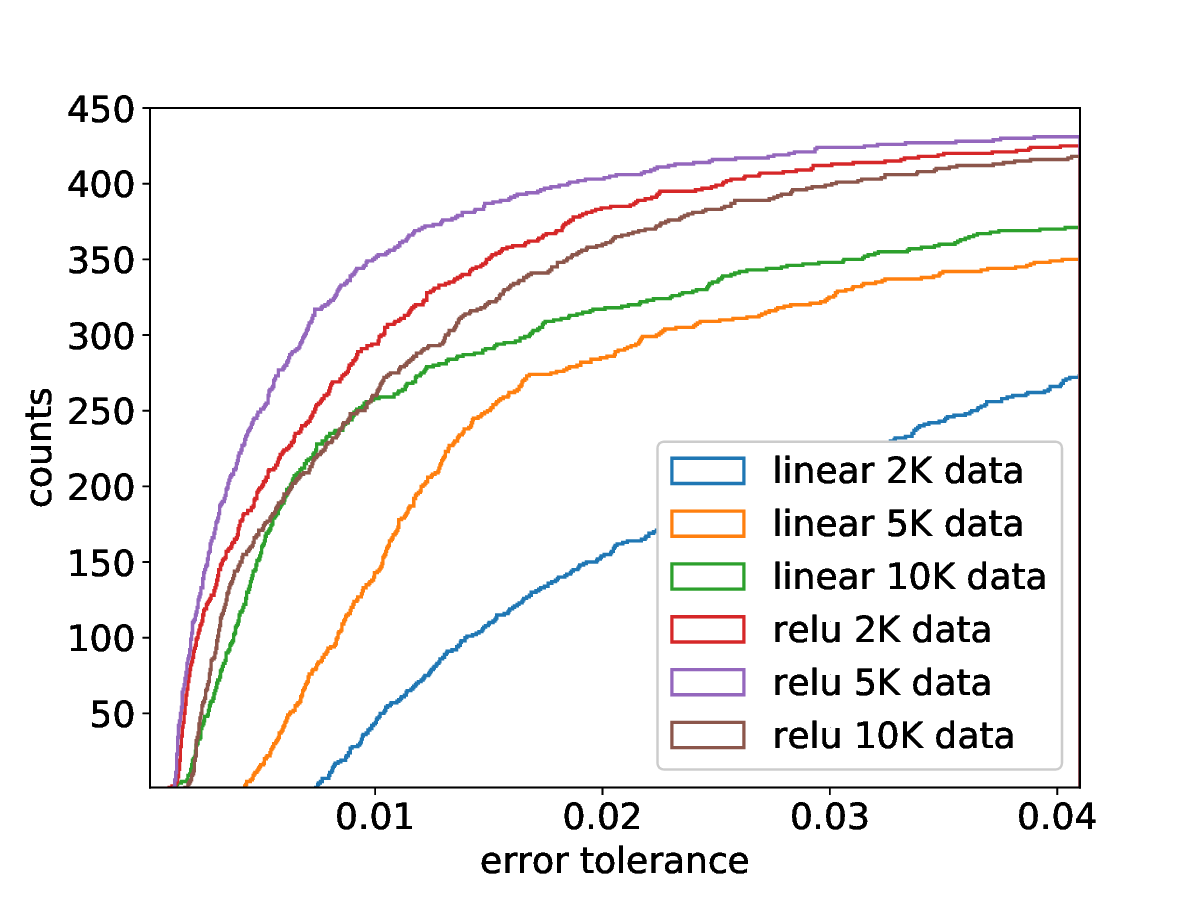}
    \caption{A comparison of the $e(\sigma, c)$ of linear and ReLU  $\mathsf{JNet(6,1)}$'s,  trained with different number of examples drawn from the same distribution, $\mathtt{D_t^p}(0.2)$. The ReLU \texttt{JNet}s exhibit better accuracy in the regime of larger wave peed and higher wave numbers.   }
    \label{fig:proxy_error_comps_dt02}
\end{figure}

We report a set of comparisons of $e(\sigma,c)$ in Figures~\ref{fig:proxy_error_comps}~and~\ref{fig:proxy_error_comps_dt02} for $\Delta t^*=0.1$ and $0.2$ respectively.
The results are obtained from the 6-level linear and ReLU \texttt{JNets}, with the parameters reported in Table~\ref{tab:disparam}.
The networks trained on $\mathtt{D_t}(0.1)$ and $\mathtt{D_t^p}(0.2)$.
In Figure~\ref{fig:proxy_error_comps},
from the two rows of error images, we observe that while linear networks produce roughly comparable errors to those from the ReLU networks in the region of lower wave speed and smaller $\sigma^{-1}$.
On the opposite regime (larger $c$ and $\sigma^{-1}$), the nonlinear ReLU networks perform better. 
As discussed in Section~\ref{sec:correction_op}, sharper wave fields and faster wave speeds result in larger nummerical dispersion errors. 
The approximation power in the nonlinear networks seems to contribute the better performance. 
We also notice that the ReLU networks trained with 10,000 examples seems to achieve slightly larger minimal errors than those by the networks trained with 5,000 examples. 
We will see in Section~\ref{sec:parareal} that a more significant number of training examples does result in better generalization capacity for the ReLU networks and perhaps not surprisingly have relatively little impact on the linear networks. 
Figure~\ref{fig:proxy_error_comps_dt02} presents a different situation.
With $\Delta t^*=0.2$, the difference between the coarse solver and the fine solver is larger. Here we see a clear advantage in using the nonlinear networks.

\begin{figure}
    \centering
    \includegraphics[width=0.8\linewidth]{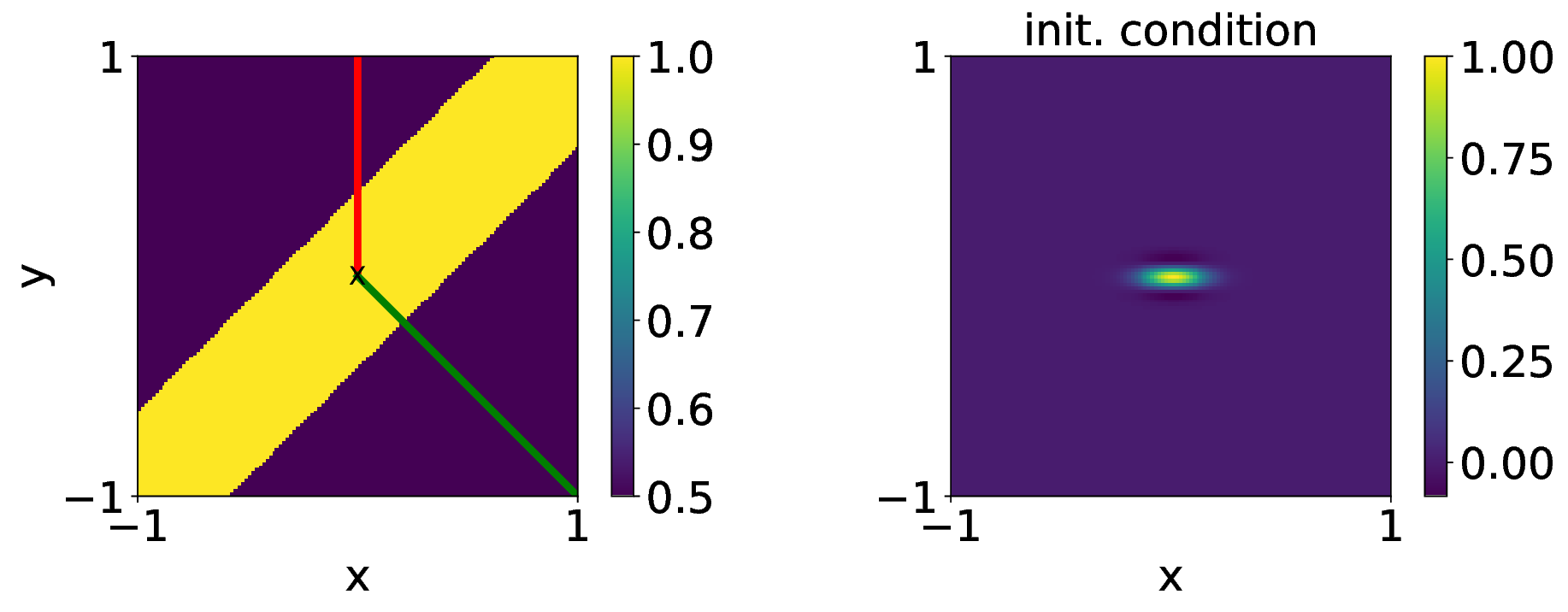}
    \caption{Left: The wave speed profile of a model layered medium. Right: The initial wave field $\mathbf{u}_0 =( \cos(8\pi y)\exp^{-25x^2-250y^2}, 0)$. }
    \label{fig:layer-test-medium-init}
\end{figure}

\subsection{Numerical dispersion and refraction}\label{sec:J3-J6-dispersion-comp}

In this section, we study the enhanced solvers in a wave medium described by piecewise constant wave speed of moderate contrast.

\begin{figure}
    \centering
    \includegraphics[width=0.98\linewidth]{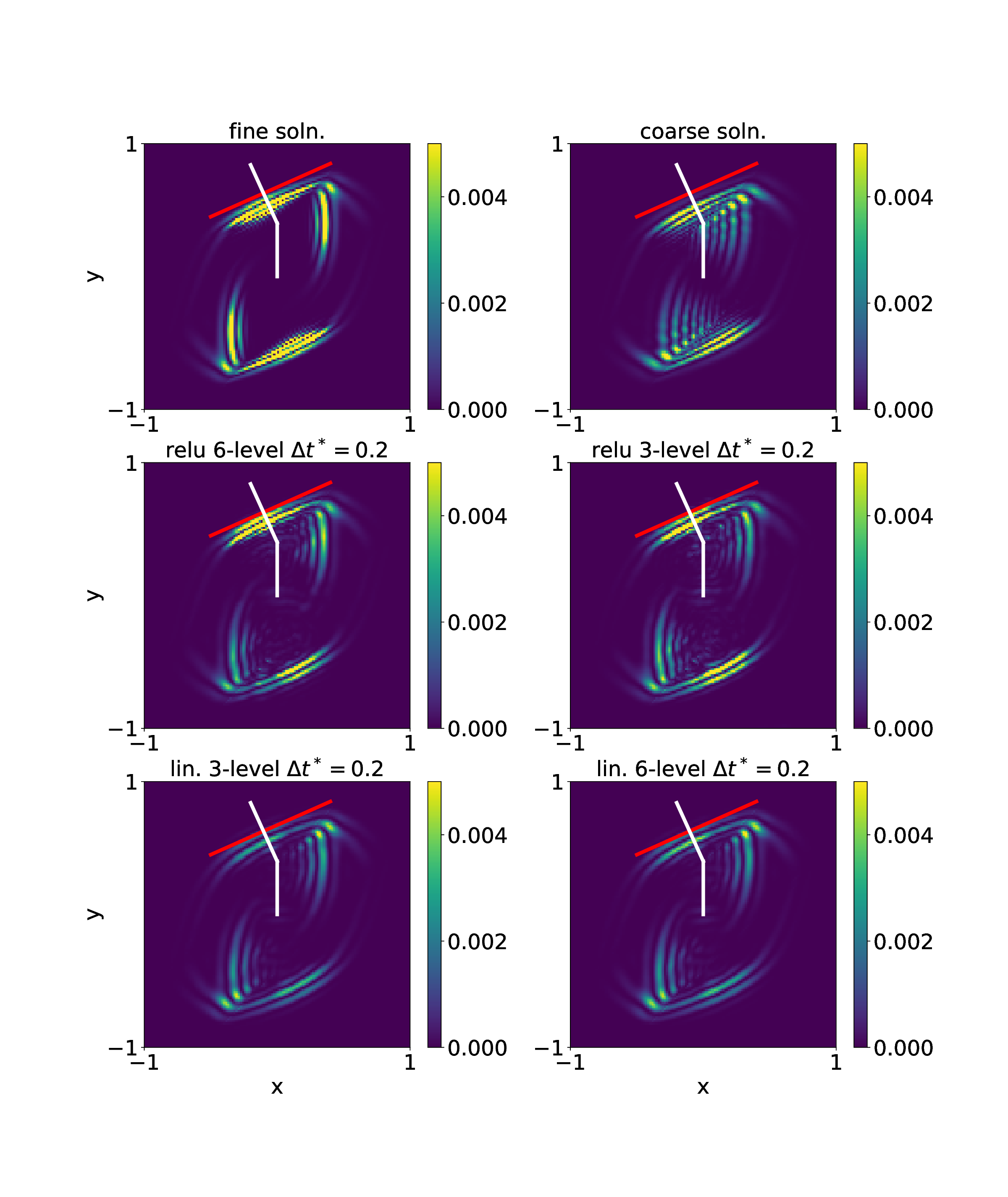}
    \caption{The proposed models correct reflected and refracted waves. }
    \label{fig:J3-J6-layerc-energy}
\end{figure}

Figure~\ref{fig:layer-test-medium-init} shows the model layered medium and the anisotropic initial Gaussian pulse. We apply different solvers to propagate that initial wave field.
The wave fronts associated with different wavenumbers travel at different speeds and transmit through the interfaces (the wave speed discontinuities) at different times. 
Figure~\ref{fig:J3-J6-layerc-energy} shows the energy of the propagated wave fields. 
The networks used in these simulations are trained on $\mathtt{D_t^P}.$
The proposed networks can provide a certain level of non-trivial correction to the dispersion and refraction in the coarse solution.  

\begin{figure}
    \centering
    \includegraphics[width=1\linewidth]{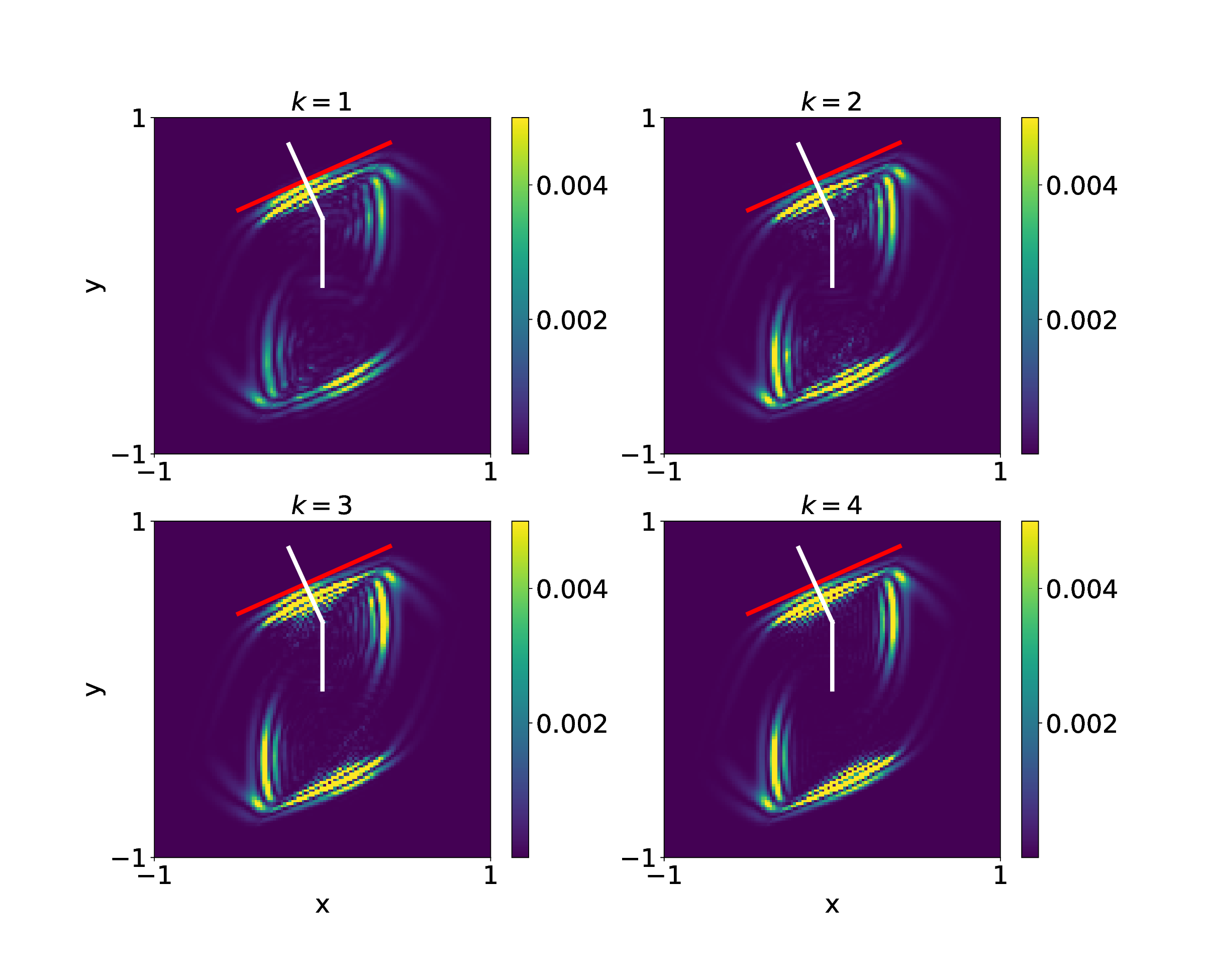}
    \caption{Wave energy field of solutions obtained from the parareal iterations, enhanced by a $\mathsf{JNet(3,1)}$. The red line helps compare the refracted wave fronts. The green line indicates the wave trajectory. The second green line segments satisfies Snell's law.}
    \label{fig:refraction-study-parareal-correction}
\end{figure}

\subsection{More complex wave speeds}
We demonstrate the effectiveness of the trained neural networks in proving the accuracy of coarse solver. 
We compute the evolution of the wave field 
$$\mathbf{u}_{n+1}= \theta_{\Delta t^*} \mathcal{G}_{\Delta t^*}\mathcal{R}\mathbf{u}_n =  \Lambda^{\dag}_{\delta x}\mathcal{H}_w(\Lambda_{\Delta x}{\mathcal{G}}_{\Delta t^*}\mathcal{R} \mathbf{u}_n,c),~~ n< N=\frac{T}{\Delta t},
$$ 
with the initial condition 
$$
    \mathbf{u}_0=(e^{-200(x^2+y^2)}, 0).
$$

\begin{figure}
    \centering
    \includegraphics[width=0.85\linewidth]{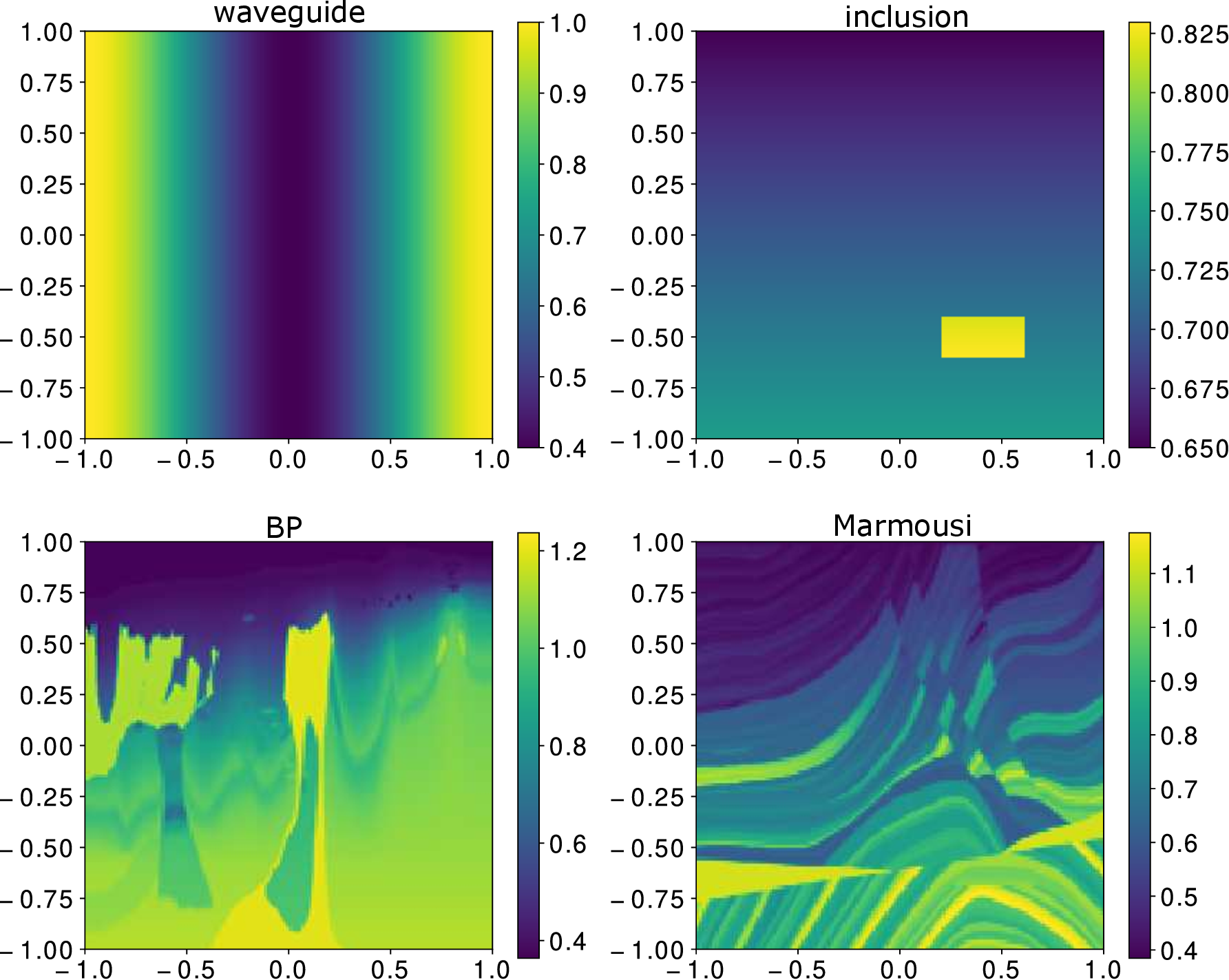}
    \caption{Wave speed models used in the experiments.}
    \label{fig:wave_speed_profiles}
\end{figure}

We consider four test media: 
\begin{itemize}
    \item Waveguide: $c(x,y) = 0.7-0.3\cos(\pi x)$,
    \item Inclusion: $c(x,y) = 0.7+0.05y+0.1\chi_{0.2<x<0.6,0.4<y<0.6}$,
    \item BP model (rescaled):  $c(x,y) = c_{BP}(x,y)/4$,
    \item Marmousi model (rescaled): $c(x,y)= c_{Marmousi}(x,y)/4$.
\end{itemize}
These four models are illustrated in Figure~\ref{fig:wave_speed_profiles}.
Note that these examples are not present in the training data set. 
The parameters used in the solvers are reported in Table~\ref{tab:disparam}.

\begin{figure}
    \centering
    \includegraphics[width=1\linewidth]{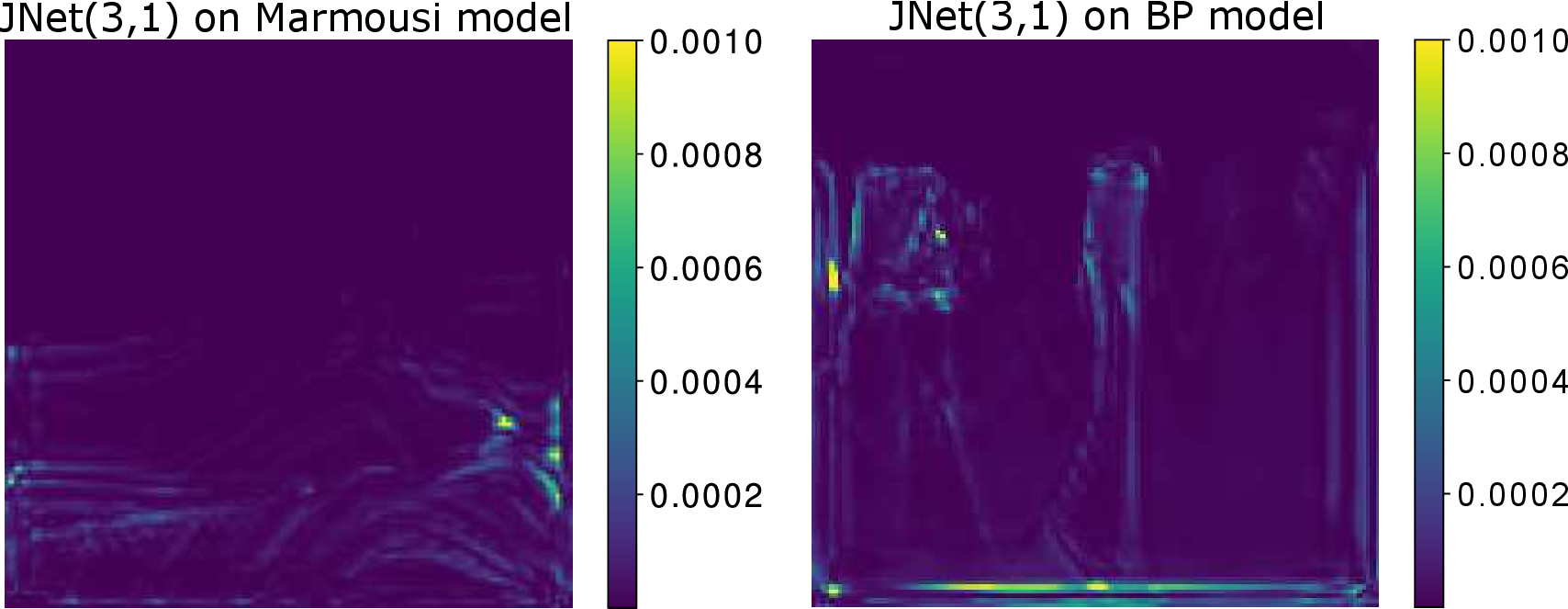}
    \caption{Phantom energy produced by the networks with zero input wave energy. Left: Marmousi model. Right: BP model. The networks are trained on $\mathtt{D_t^p(0.2)}$.}
    \label{fig:phantom_energy}
\end{figure}

In Figure~\ref{fig:phantom_energy} we show the ``phantom" energy in the output of $\mathcal{H}_w\mathbf{0}.$ While this can be regarded as an undesirable effect, it does improve the amplitude of the propagated wave field. We also observe that the trained 6-level network produces drastically reduced amount of ``phantom" energy.

\section{Improvement by parareal iterations}\label{sec:parareal}

One of the shortcomings of employing a current deep learning strategy for scientific computing tasks is the lack of rigorous error estimates. 
Our idea is to use the proposed deep learning strategy in a self-improving loop, in which one can show that the computed errors will contract along the iterations. Here, we adopt the parareal framework, originally proposed in \cite{parareal-LMT01}, as the potentially self-improving looping mechanism. 

Parareal methods can be regarded as a fixed point iteration that involves the coarse and the fine solvers:
\begin{eqnarray}
    \mathbf{u}^{k+1}_{n+1} &=& \tilde{\mathcal{G}}_{\Delta t^*} \mathbf{u}^{k+1}_{n} + \mathcal{F}_{\Delta t^*} \mathbf{u}^{k}_{n} - \tilde{\mathcal{G}}_{\Delta t^*} \mathbf{u}^{k}_{n},~k=0,1,2,\dots, 
\end{eqnarray}
where the loop is enumerated in $k$ and the subscript $n$ denotes the time stepping.  The initial conditions for the iterations are defined by
\begin{eqnarray}\label{eq:plain-parareal-init}
    \mathbf{u}^{0}_{n+1} &=& \tilde{\mathcal{G}}_{\Delta t^*} \mathbf{u}^{0}_{n},~~~~~~n=0,1,\dots,N-1.
\end{eqnarray}
In the context of this paper,  
\[\mathbf{u}^n_k :=  (u^n_k,\partial_t u^n_k)\]
is defined on the fine grid, and 
$\tilde{\mathcal{G}}_{\Delta t^*}:=\mathcal{I}\mathcal{G}_{\Delta t^*}\mathcal{R}.$

The convergence of parareal methods are analyzed for different contexts in several papers. See e.g. \cite{ANT-thetaparareal},\cite{Bal2005},\cite{gander2007analysis},\cite{mikioInfluence18},\cite{Ruprecht18}.
We adhere to a simple version, which states that at the discrete level with the vector $\infty$-norm, the error amplification factor for the parareal iteration is the product of  $||\mathcal{F}_{\Delta t^*}-\tilde{\mathcal{G}}_{\Delta t^*}||$ and $\frac{1-r^N}{1-r},$ with $r\equiv ||\tilde{\mathcal{G}}_{\Delta t^*}||$.
The iteration is convergent when the coarse solver is sufficiently accurate, depending on the magnitude of $||\tilde{\mathcal{G}}_{\Delta t^*}||$ and the size of the iterative system, $N$.
It is a challenge to construct a stable and convergent parareal method for purely hyperbolic problems, especially when the coarse solver is under-resolving the medium. 

As we have seen in the serial computation examples in Section~\ref{sec:netexample}, our network reduces the errors in the coarse solutions by 90\% in very challenging media. 
Hence, we expect the parareal iteration involving the enhanced coarse solver to be more stable and, in turn, be used to improve approximation accuracy.



Following the framework proposed as the $\theta$-parareal schemes \cite{ANT-thetaparareal}\cite{nguyen2020stable}, the focus of this section rest on schemes of the form: 
\begin{equation} \label{eq:nnparareal}
    \mathbf{u}^{k+1}_{n+1} = \theta_{\Delta t^*}({\mathcal{G}}_{\Delta t^*} \mathcal{R}\mathbf{u}^{k+1}_{n}, c) + \mathcal{F}_{\Delta t^*} \mathbf{u}^{k}_{n} - \theta_{\Delta t^*} ({\mathcal{G}}_{\Delta t^*} \mathcal{R}\mathbf{u}^{k}_{n},c),
\end{equation}
where the operator $\theta_{\Delta t^*}$ under consideration include
\begin{itemize}
    \item $\theta_{\Delta t^*}$ is defined by a \texttt{JNet}, trained offline. \LL{This is the main focus of this paper.}
    \item $\theta_{\Delta t^*}\equiv\Lambda^\dag\Omega$ where $\Omega$ is a unitary matrix derived from the online data $\{\mathbf{u}^k_n\}.$ \LL{This approach is used to generate training database $\mathtt{D_t^p}$.}
\end{itemize}

\subsubsection*{The neural network approach}

The central point of this paper is 
\[
\theta_{\Delta t^*}(\mathcal{G}_{\Delta t^*}\mathbf{u}, c)\equiv \Lambda^\dag\mathcal{H}_{w}({\Lambda\mathcal{G}}_{\Delta t^*} \mathcal{R}\mathbf{u}^{k}_{n-1}, c),
\] for $\mathcal{H}_{w}$ is a trained neural network discussed above.
We initialize the parareal iterations by the enhanced coarse solutions: 
\begin{equation} \label{eq:nnparareal-init-enhanced}
    \mathbf{u}^{0}_{n+1} = \theta_{\Delta t^*}({\mathcal{G}}_{\Delta t^*} \mathcal{R}\mathbf{u}^{0}_{n}, c),~~~n=0,1,\dots,N=T/{\Delta t^*-1.}
\end{equation}

\subsubsection*{The Procrustes parareal method}


The Procrustes parareal method is defined as follows:
\begin{equation} \label{eq:Procrustes-parareal}
    \mathbf{u}^{k+1}_{n+1} = \Lambda^\dag\Omega^k\tilde{\mathcal{G}}_{\Delta t^*} \mathbf{u}^{k+1}_{n} + \mathcal{F}_{\Delta t} \mathbf{u}^{k}_{n} - \Lambda^\dag\Omega^k \tilde{\mathcal{G}}_{\Delta t^*}\mathbf{u}^{k}_{n},~~k=0,1,\dots,
\end{equation}
where $\tilde{\mathcal{G}}_{\Delta t^*}:=\mathcal{I}\mathcal{G}_{\Delta t} \mathcal{R},$
and for $k>0$, $\Omega^k$ solves the following optimization problem:
\begin{equation}\label{eq:Procrustes-optimization}
    \min_{\Omega\in\mathbb{R}^{N\times N}} ||\Omega \mathsf{G}_k -\mathsf{F}_k||_F^2~~~s.t.~~~ \Omega^T\Omega=I.
\end{equation}
The data matrices $\mathsf{G_k}$  and $\mathsf{F}_k$ are defined as follows:
\begin{itemize}
    \item The $n$-th column of $\mathsf{G_k}$ corresponds to the vector $\Lambda\tilde{\mathcal{G}}_{\Delta t^*} \mathbf{u}^k_n$, $n=0,1,\dots,N.$
    \item The $n$-th column of $\mathsf{F_k}$ corresponds to the vector $\Lambda\mathcal{F}_{\Delta t} \mathbf{u}^k_n$, $n=0,1,\dots,N.$
\end{itemize}
For $k=0$, $\mathbf{u}^0_n$ are defined in \eqref{eq:plain-parareal-init}.
We call this method the Procrustes parareal method since
 $\Omega^k$ is constructed by solving the Orthogonal Procrustes Problem \cite{golub2012matrix}.

The Procrustes parareal algorithm is introduced in \cite{nguyen2020stable}. In the original version, the fine and coarse solvers communicate on the coarse grid. This means that the data matrices records the wave fields on the coarse grid and $\Omega^k$ is smaller. In contrast, in the version presented above, the the solvers communicate on the fine grid. This allows the resulting algorithm to correct for details corresponding to higher wave numbers.
\LL{One disadvantage of this method is the potential scalability issue due to $\Omega$ being large and dense for large computational domains. However, it may be combined a domain decomposition approach (in space) to reduce the problem size.}
 
Both versions of the Procrustes parareal method are remarkably stable for wave propagation. We point out again, that due to this stability, we are able to use this method to generate training examples that lead to more performant neural networks. See Section~\ref{sec:netexample}.

\subsection{Numerical simulations and comparisons}
First, we follow up the dispersion study in Section~\ref{sec:J3-J6-dispersion-comp} and 
present the wave energy field computed in the first four parareal iterations in Figure~\ref{fig:refraction-study-parareal-correction}. There we observe the effect of the additive correction of the parareal updates --- with accurate phase, the amplitudes are restored effectively.

\begin{figure}
    \centering
    \includegraphics[width=1\linewidth]{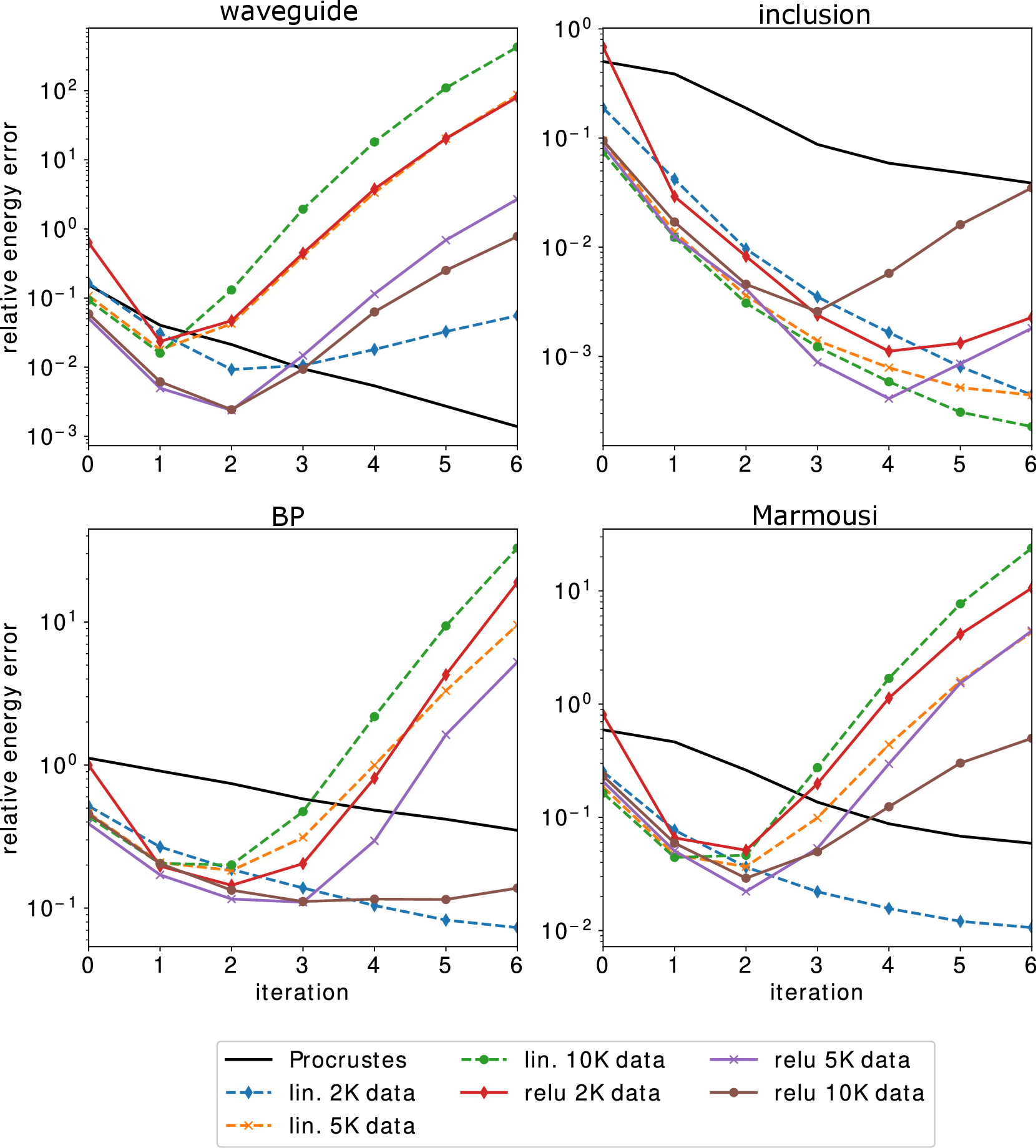}
    \caption{A comparisons of parareal iterations applied to the \texttt{JNets} trained on data sets of different sizes. The training sets are uniformly chosen subsamples of $\mathtt{D_t}(0.1).$ The errors are evaluated at $T=1.5$. }\label{fig:comprehensive-comp}
\end{figure}

\begin{figure}
    \centering
    \includegraphics[width=1\linewidth]{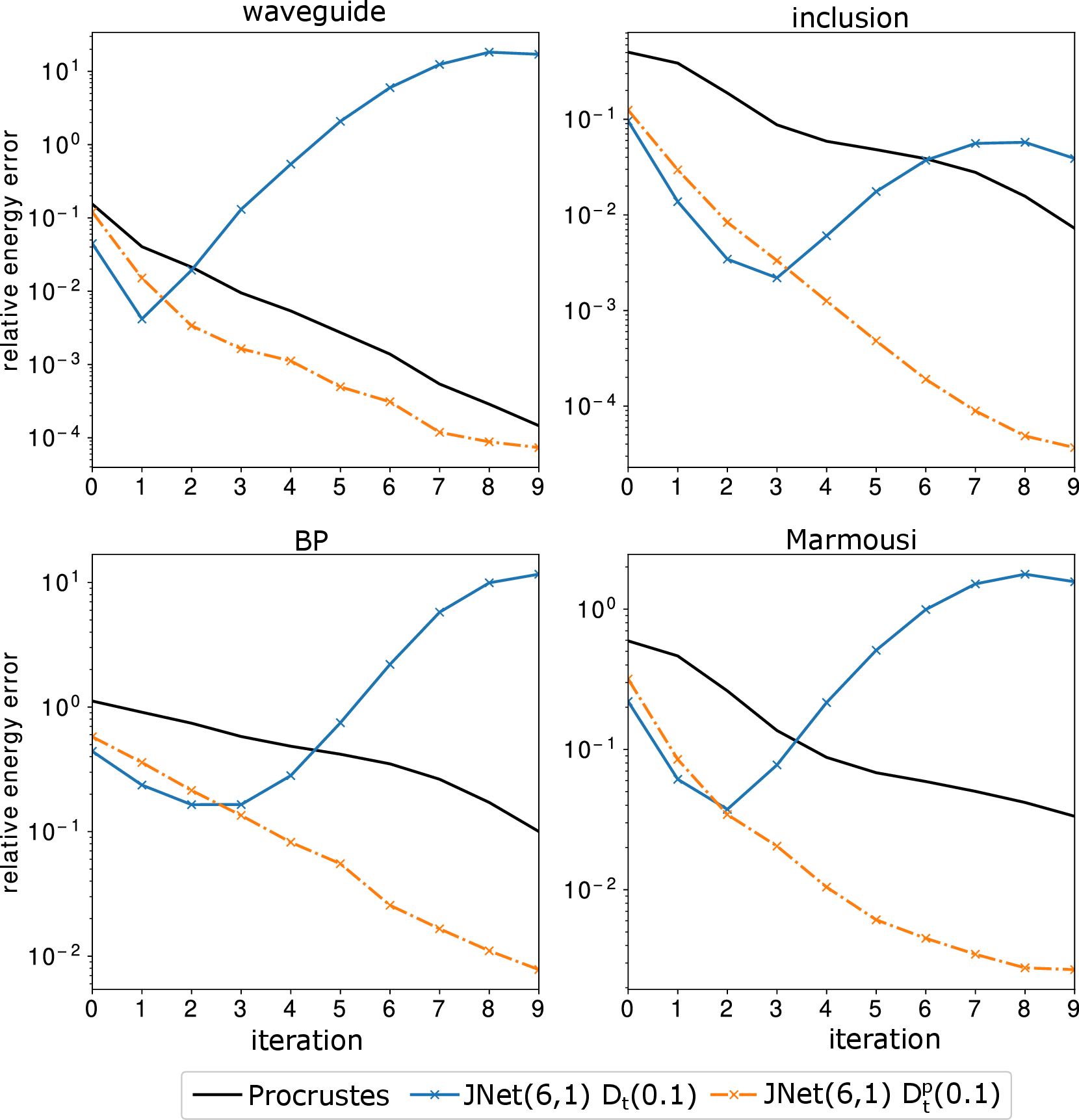}
        \caption{Effects of the training data sets. The errors are evaluated at $T=1.5$. }
    \label{fig:J6-comp-of-training-datasets}
\end{figure}

\begin{figure}
    \centering
    \includegraphics[width=1\linewidth]{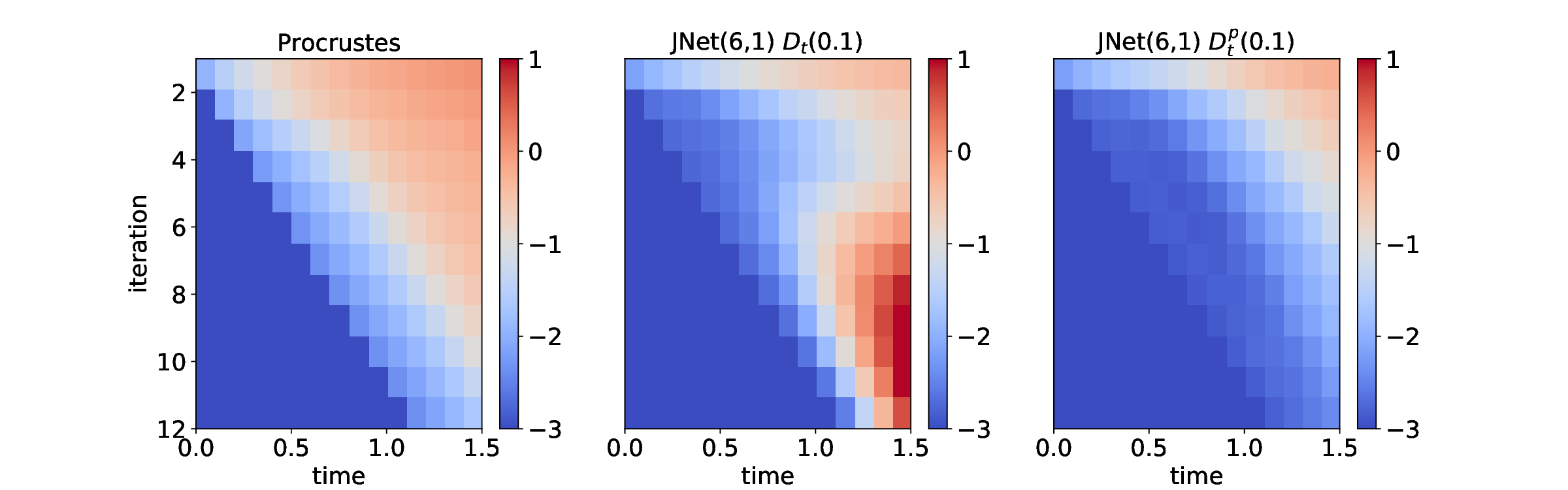}
    \caption{\LL{The effect of training sets in the computed relative energy errors, at different times of propagation and different parareal iterations. The simulations are performed on the BP model.}}
    \label{fig:J6-comp-heatmap}
\end{figure}

\begin{figure}
    \centering
    \includegraphics[width=1\linewidth]{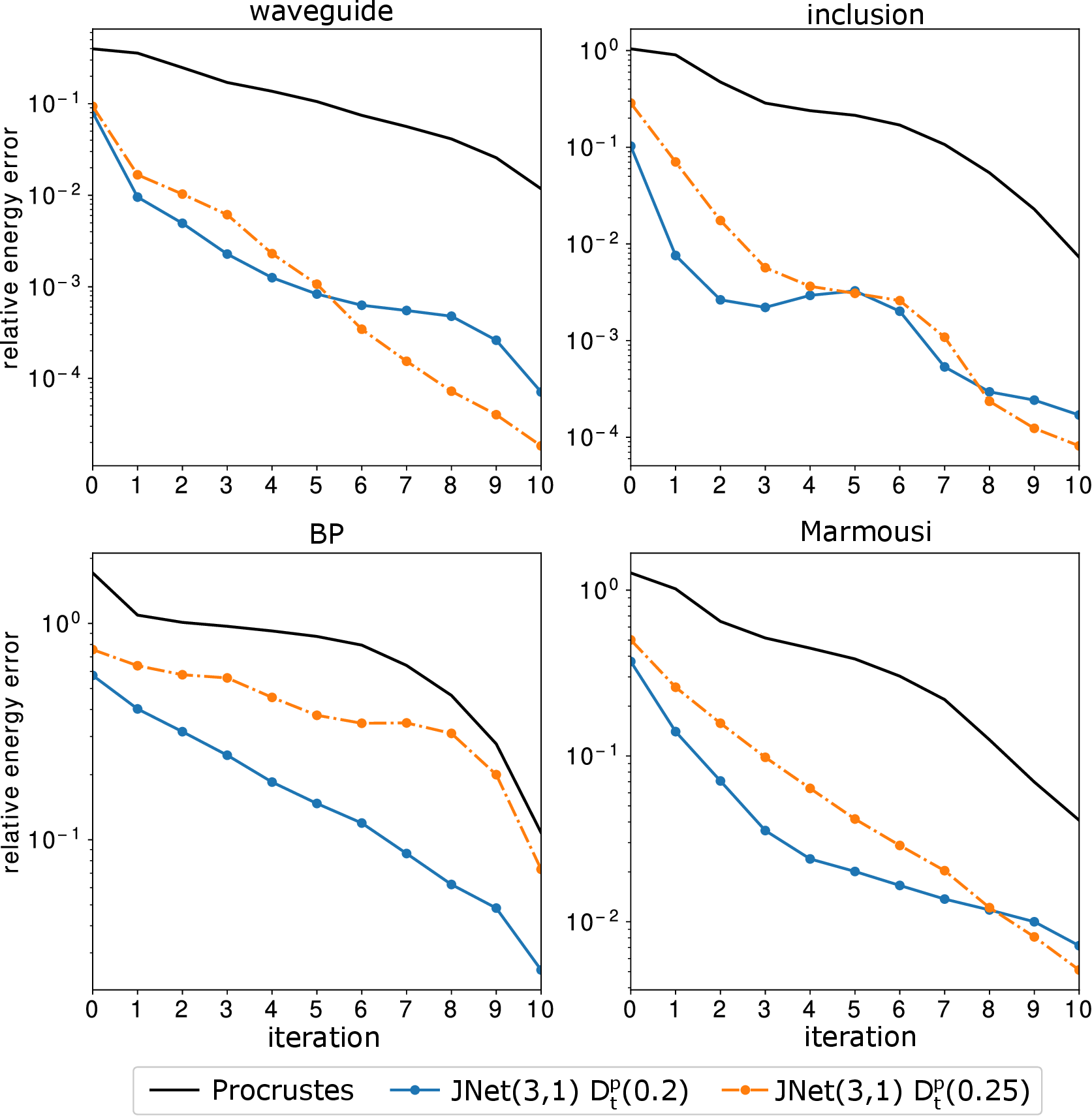}
    \caption{Effect of the size of $\Delta t^*$. Accuracy and convergence of the parareal iterations, using 3-level \texttt{JNet}s trained on $\mathtt{D_t^p}(0.2)$ and $\mathtt{D_t^p}(0.25)$. 
    The final time of propagation is $T=3.0$. }\label{fig:J3-fixed-T-different-dt*}
\end{figure}

\begin{figure}
    \centering
    \includegraphics[width=1\linewidth]{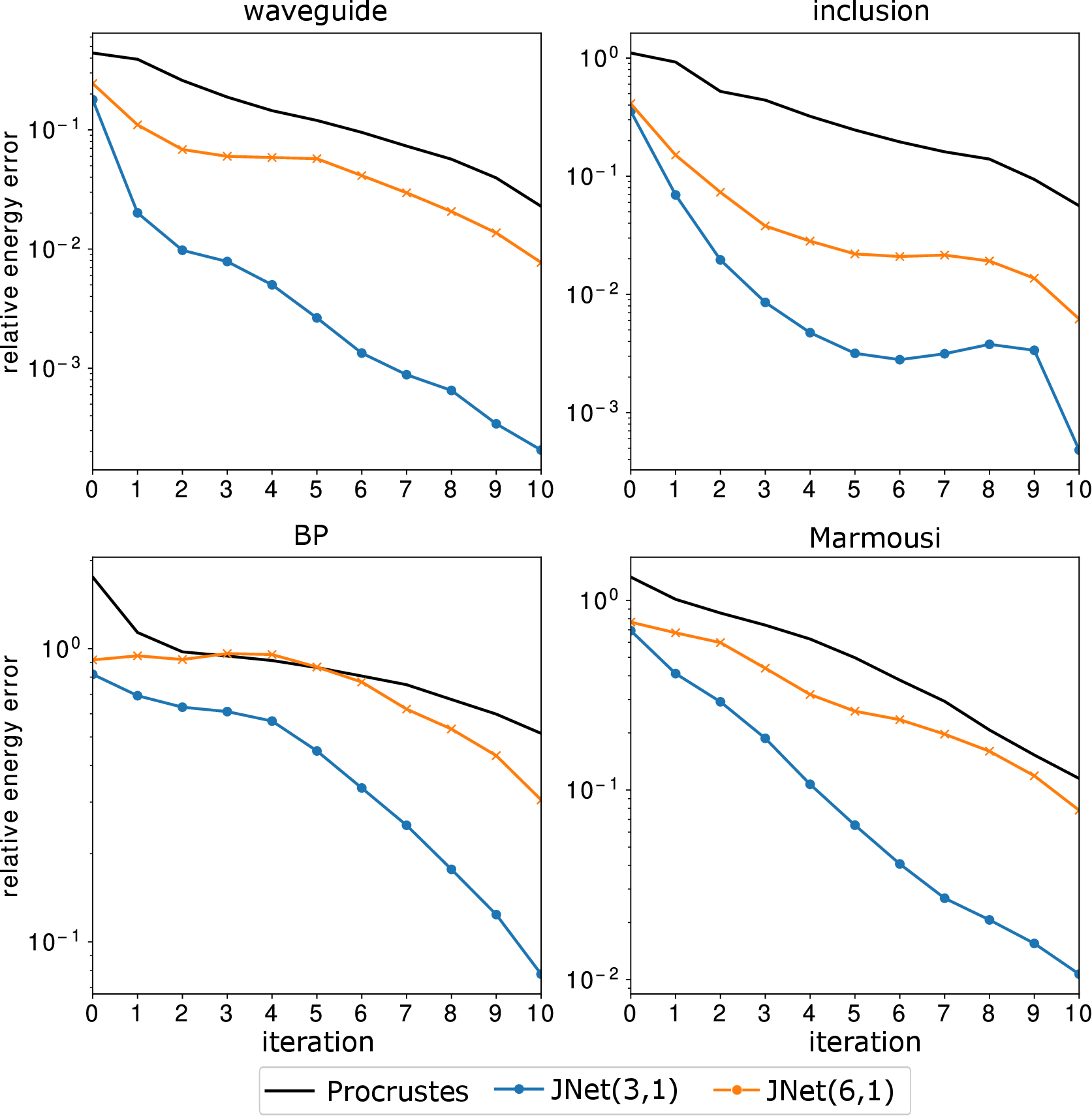}
    \caption{A comparison of parareal iterations enhanced by $\mathsf{JNet(3,1)}$ and $\mathsf{JNet(6,1)}$.  The networks are trained on $\mathtt{D_t^p}(0.2).$ The final time of propagation is $T=3.2$.}
    \label{fig:parareal-J3-J6_comp_fixed-T}
\end{figure}

\begin{figure}
    \centering
    \includegraphics[width=1\linewidth]{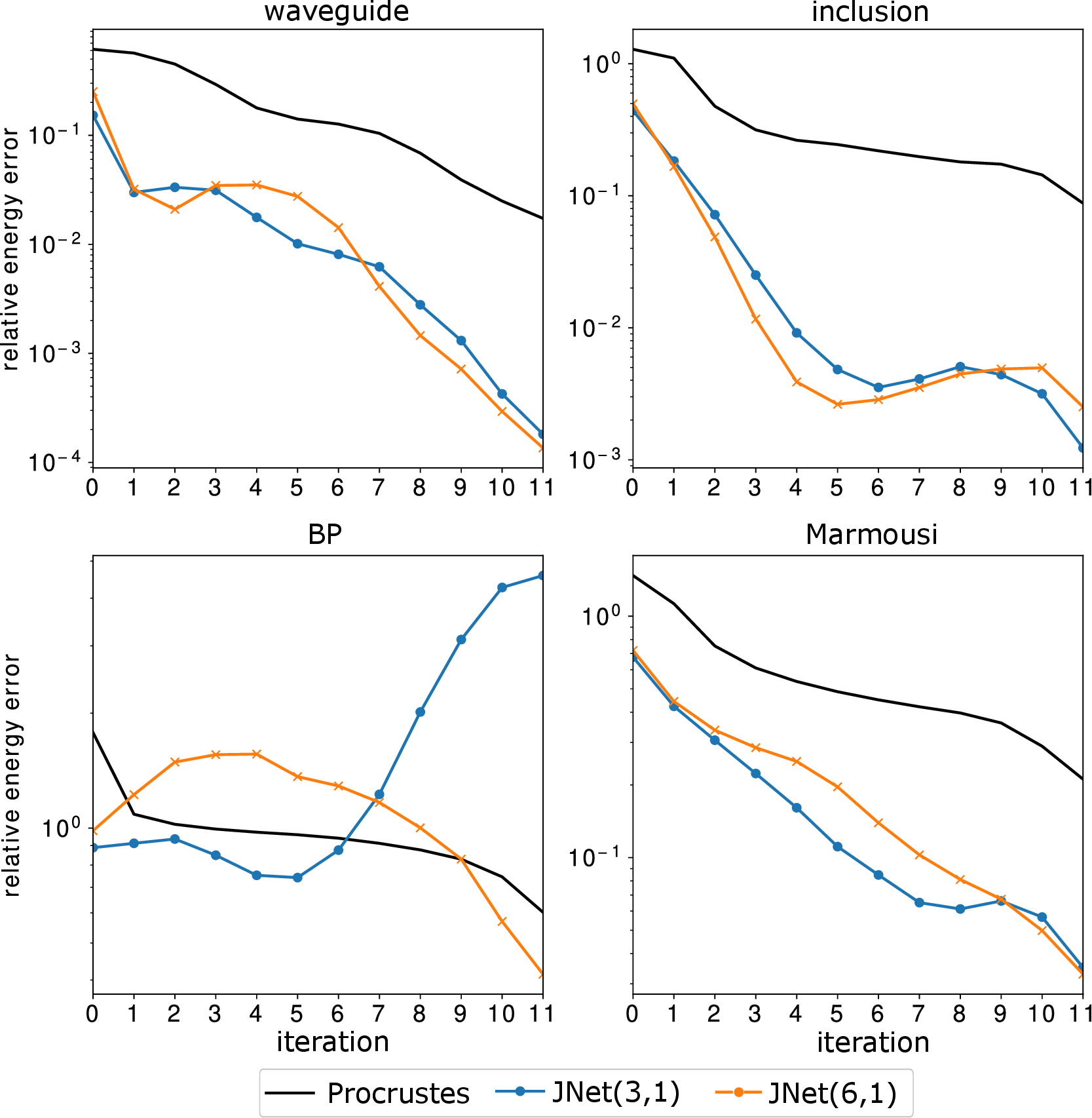}
    \caption{Accuracy and convergence of the parareal iterations involving networks of different depths. The 3-level and 6-level \texttt{JNet}s are trained on $\mathtt{D_t^p}(0.25)$. The final propagation time is
    $T=4.0$.}\label{fig:comp-Procrustes-vs-J6Net-dt=0.25}
\end{figure}

In Figures \ref{fig:comprehensive-comp}-\ref{fig:comp-Procrustes-vs-J6Net-dt=0.25}, we observe that the initial errors ($k=0$) 
dominated by the effects of the restriction operator in $\mathcal{G}_{\Delta t^*}\mathcal{R}$, being repeatedly applied in the time stepping process. However, as $k$ increases, the effect of having better phase errors in the \texttt{Jnet}s can be observed. In some situations, the advantage of having more dissipation can also be observed.  

\paragraph{Linear vs. nonlinear networks.} We present a comprehensive accuracy study in Figure~\ref{fig:comprehensive-comp}, comparing the solvers that we have discussed so far.
There, we observe the importance of the training data size for the nonlinear networks.
We also notice that the linear networks perform surprisingly well in the studied setup (particularly the linear network trained with on 2000 examples).
One possible explanation is that the operator to be approximated, 
\[
\mathcal{F}_{\Delta t^*}\mathcal{I}\mathcal{G}_{\Delta t^*}^{-1},
\]
is closer to the identity matrix 
for $\Delta t^*=0.1$ and the class of $c$.
Ergo, linear regression may yield reasonable approximations. Indeed, we have seen in Figure~\ref{fig:proxy_error_comps_dt02} that linear networks are inferior for larger $\Delta t^*$.

\paragraph{The effects of the training sets: $\mathtt{D_t}$  vs. $\mathtt{D_t^p}$.}
Figure~\ref{fig:J6-comp-of-training-datasets} shows a comparison of the networks trained on $\mathtt{D_t}$ and $\mathtt{D_t^p}.$
We observe that the network trained on $\mathtt{D_t^p}$ yields better performance after a few iterations. 
\LL{We present a more detailed comparison of the simulations on the BP model in Figure~\ref{fig:J6-comp-heatmap}.}
We speculate that the wave fields in $\mathtt{D_t^p}$ are
closer to the ones encountered in the simulations performed in this experiment --- the wave fields are computed by parareal style updates. \LL{In other words, after 5 parareal iterations, the comptued wave fields may become too different from the training examples collected in $\mathtt{D_t}$.} \LL{The deterioration of the errors may thus be related to the stability of a trained neural network for making ``out-of-distribution" inferences, see \cite{HTW2022}.}

\begin{figure}
    \centering
    \includegraphics[width=1\linewidth]{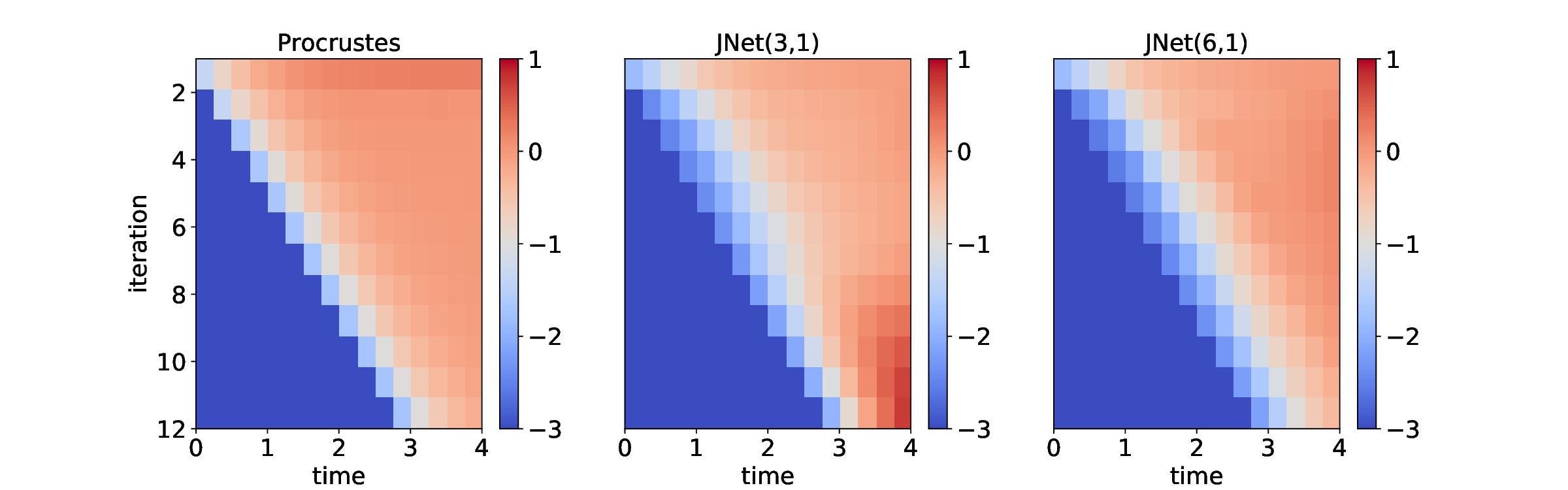}
    \caption{Relative energy error at different times of propagation and different parareal iterations on the BP model. The parareal iterations involves \texttt{JNet}s networks trained on $\mathtt{D_t^p}(0.25)$.}
    \label{fig:comp-Procrustes-vs-J6Net-dt=0.25_heatmap}
\end{figure}

\paragraph{The effect of $\Delta t^*$ to the accuracy.}
In Figure~\ref{fig:J3-fixed-T-different-dt*} we compare the accuracy of the 3-level \texttt{JNet}s trained for $\Delta t^*=0.2$ and $0.25.$ 
We first observe from the errors at $k=0$ that the enhanced solver with the smaller $\Delta t^*$ is more accurate. 
Except for the BP-model, the performance of the parareal iterations with the two networks seem comparable.
\LL{We continue to study the simulations related to the BP model in Figure~\ref{fig:comp-Procrustes-vs-J6Net-dt=0.25_heatmap}, where we see the larger errors in $\mathsf{JNet}(3,1)$ (trained with $\mathtt{D_t^p}(0.25)$) accumulate and take a more pronounced effect at a later time.}

\paragraph{$\mathsf{JNet}(3,1)$ vs. $\mathsf{JNet}(6,1)$.} In Figure~\ref{fig:parareal-J3-J6_comp_fixed-T}, we present a comparison of the performance of the proposed parareal correction, \eqref{eq:nnparareal}, using the three-level and the six-level networks, both trained on $\mathtt{D_t^p}$. The iterations with the three-level network appear to yield smaller errors because the three-level network has larger numerical dissipation (See Figure~\ref{fig:J3-J6-layerc-energy} and Section~\ref{sec:J3-J6-dispersion-comp}). 

In Figure~\ref{fig:comp-Procrustes-vs-J6Net-dt=0.25}, we compare the Procrustes approach defined in \eqref{eq:Procrustes-parareal} and the $\mathsf{JNet}(3,1)$ and $\mathsf{JNet}(6,1)$-enhanced coarse solvers, for a larger $\Delta t^*=0.25$ and larger final time $T=4.0$. We notice that, for the BP model, the final time error computed by the  $\mathsf{JNet}(3,1)$-enhanced solver started to increase after the fifth parareal iteration. 
\LL{Comparing to  Figure~\ref{fig:parareal-J3-J6_comp_fixed-T}, the difference is the larger $\Delta t^*$ and final time. However, Figure~\ref{fig:comp-Procrustes-vs-J6Net-dt=0.25_heatmap} reveals that the simulation involving $\mathsf{JNet}(3,1)$ with $\Delta t^*=0.25$ already deteriorated at $t\approx 3$. This suggests that the trained $\mathsf{JNet}(3,1)$ is not sufficiently accurate for such a large step size.}

\paragraph{Further remarks.} So what can we conclude from these examples? First, the phase accuracy is of upmost importance in the performance of the parareal-style iterations. Second, the parareal scheme itself is sensitive to the amount of dissipation in the system. For systems with little or no dissipation, such as the hyperbolic system that we consider in this paper, the size of the iterative system, $N=T/{\Delta t^*}$ will have a critical role: the larger it is, the smaller the difference between $\mathcal{F}_{\Delta t^*}$ and $\mathcal{G}_{\Delta t^*}$ has to be in order to have stability and achieve the desired correction. Of course, the stability and convergence for such type of iterations is well understood, e.g. \cite{xu2017algebraic}. 
Our examples not only confirm such theories, but also demonstrate the feasibility of adopting a data-driven approach, e.g. the Procrustes approach \cite{nguyen2020stable} using online data or the deep learning approach using offline data.

Finally, we noticed that the BP model presents more difficulty for the proposed strategy. One possible cause may be the higher contrast in the BP model --- particularly the large discontinuity across the periodically extended top and bottom boundaries --- and that the training data should include more instances of wave propagating through such type of interfaces. This is a subject of future investigation.

\subsubsection{The unscaled Marmousi velocity model}\label{sec:unscaled-Marmousi-sim}
\begin{figure}
    \centering
    \includegraphics[width=0.58\linewidth]{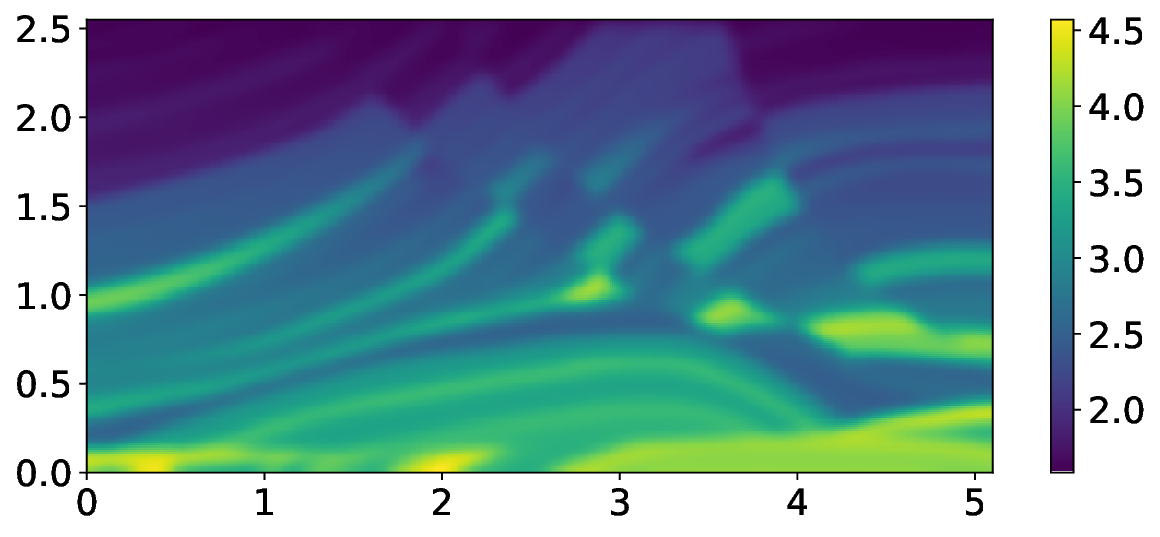}
    \includegraphics[width=1\linewidth]{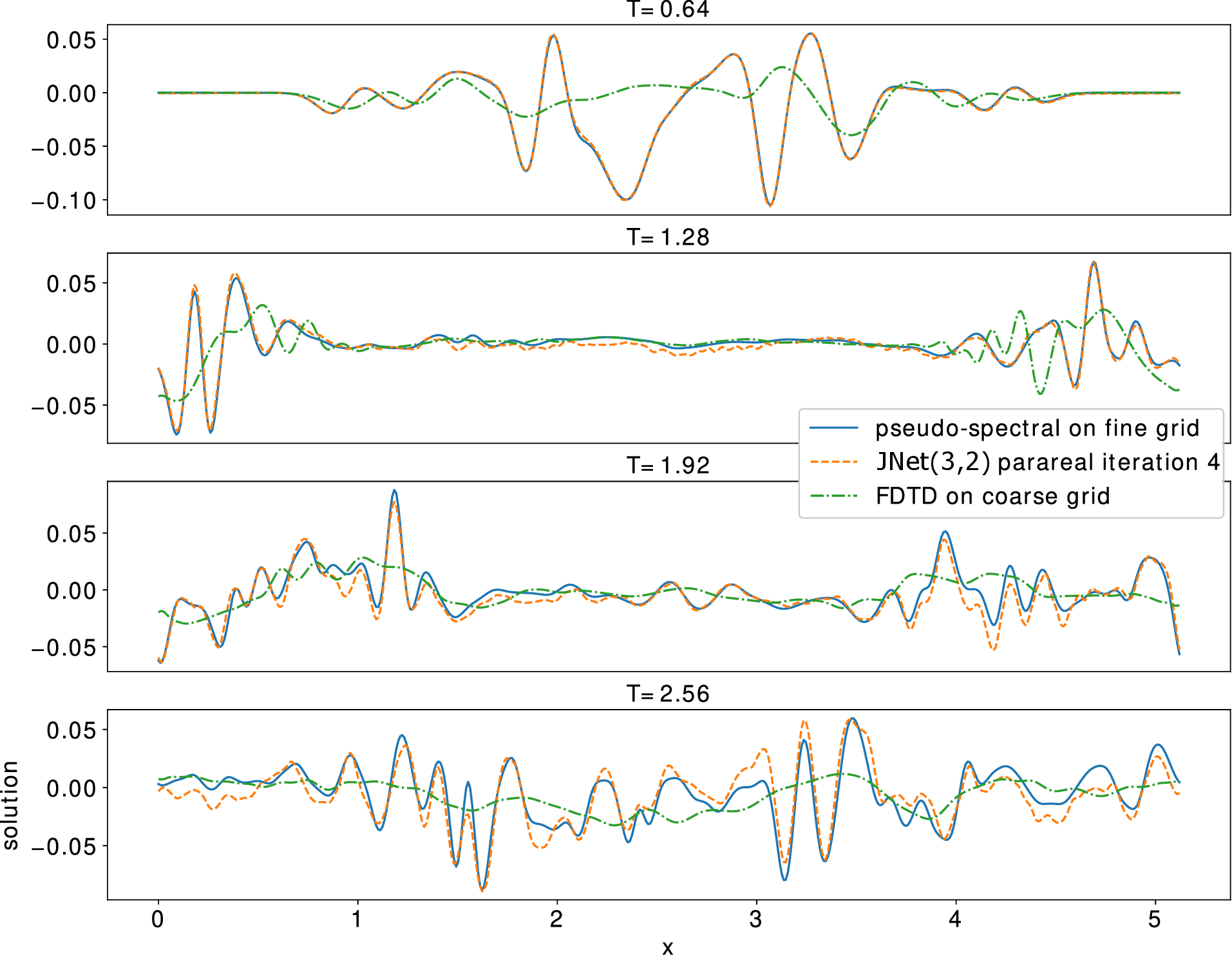}
    
    \caption{\LL{Parareal solutions enhanced by a $\mathsf{JNet}(3,2)$ on a unscaled subdomain of the mollified Marmousi model. The units are in $km$ and $second$.}}
    \label{fig:unscaled_Marm}
\end{figure}

\begin{figure}
    \centering
    \includegraphics[width=1.05\linewidth]{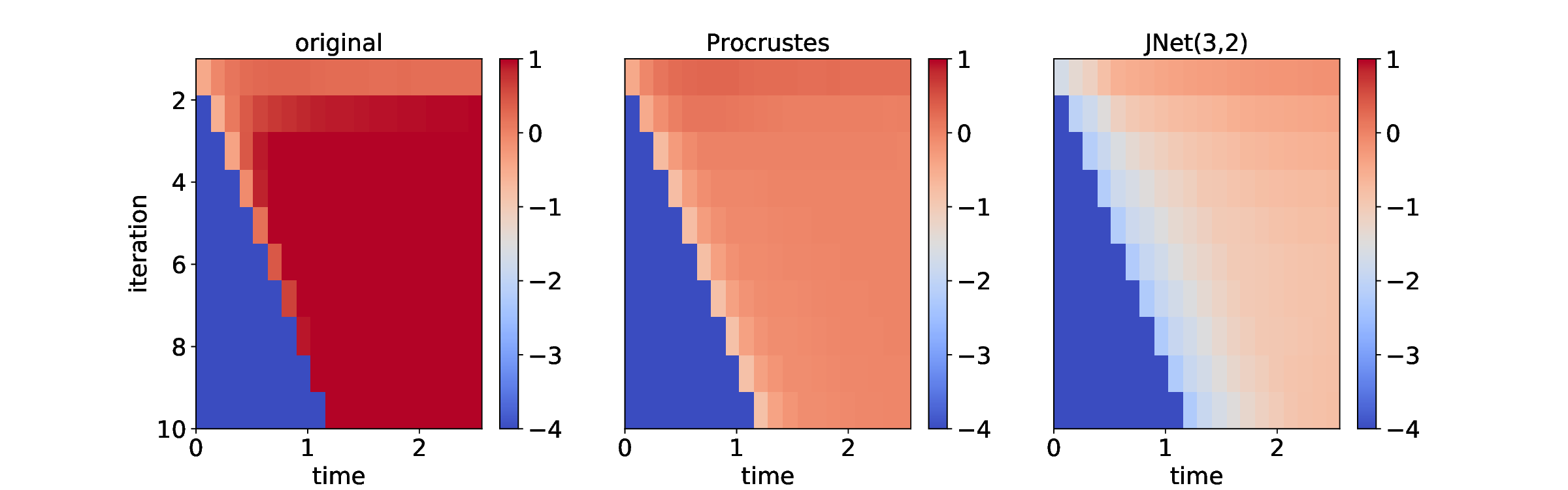}
    \caption{Parareal performance. The heatmap color shows the logarithmic scale (base $10$) of the relative energy error at different times and parareal iterations.}
    \label{fig:parareal-convergence-Marmousi}
\end{figure}

\LL{The purpose of this Section is to test the propose Deep Learning + parareal approach in a slightly more complicated setup to reveal its potential for scaling up.
In particular, the fine solver is the pseudo-spectral method and the coarse solver is the central scheme discussed in Section~\ref{sec:numerical-solver}. 
 The coarse solver uses $\Delta x = 0.04~km, \Delta t=0.004~s$. The fine solver uses $\delta x=0.01~km, \delta t=5\cdot 10^{-4}~s$. 
 }
 
\LL{We train and test the performance of $\mathsf{JNet}(3,2).$ This means that the networks 
take as input  $64\times 64\times 3$ tensors, and output a $256\times 256\times 3$. }

\LL{The networks are trained on $\mathtt{D_t^p}(0.128~s),$ which is generated by the procedure presented in Section~\ref{sec:netexample} from random crops of the Marmousi and the BP models, and initial wave field defined by the Ricker wavelets. Following a common practice, the unscaled Marmousi and BP wave speeds are mollified with a Gaussian filter with standard deviation of size $4\delta x~ (0.04km).$ The Ricker wavelets take the form
\[ 
    u_0 (x,y) \varpropto \Big(1-2\dfrac{(x-x_0)^2+(y-y_0)^2}{\sigma^2}\Big)\exp{\Big(-\dfrac{(x-x_0)^2+(y-y_0)^2}{\sigma^2}\Big)}
\]
for $\sigma^{-2} \sim \mathcal{N}(250,10)$ and have randomly selected center $(x_0,y_0)$. 
}

\LL{One of the benefits of the convolution network is the flexibility to take larger input tensors and processed them with the optimized convolution kernels, provided that the input tensors have the same physical units.}

\LL{We apply $\mathsf{JNet}(3,2)$ to enhance 
coarse wave propagation over the Marmousi model on a $512\times 256$ grid, starting from a Ricker pulse, of width $0.05~km$, locates at $(1.28~km,2.56~km)$ 
This entire wave speed is not present in the training data set.
The final time of propagation for the parareal iterations is $T=2.56~s$. }
 
\LL{In Figure~\ref{fig:unscaled_Marm}, we show the solution at the top of the computational domain ($z=2.5$) at different times. 
There we observe that the neural network greatly improves the accuracy of the coarse computation.
Figure~\ref{fig:parareal-convergence-Marmousi} shows the convergences of several parareal schemes. The original one is not stable for solving hyperbolic equations. The Procrustes approach is stable, but its error stagnates after only two iterations. We think this stagnation is due to the non-local correction operator. The $\mathsf{JNet}(3,2)$-enhanced parareal scheme significantly decreases the error in the wave energy after the first five iterations; see Figure~\ref{fig:parareal-convergence-Marmousi}. }

\LL{We recorded the average time to 
to advance a given wave field from $t_n$ to $t_n+\Delta t^*$, with $\Delta t^*=0.128 s$. On a CPU of a moderate personal workstation,  the pseudo-spectral fine solver took $4.53s$.
The feedforward of $\mathsf{JNet}(3,2)$ took $2.74s$. 
Furthermore, discounting the parallel overhead, a 
$\mathsf{JNet}(3,2)$-enhanced parareal iteration
on a 20-core environment would provide a hypothetical $33\times$ speed-up of simulation time on the wall-clock.
We believe better implementations and utilizing special hardwards, such as GPUs or TPUs as in \cite{kochkov2021machine}, will significantly reduce the computing time.}


\section{Summary}
In this paper, we presented a deep learning approach for improving the accuracy of wave propagation by an under-resolving numerical solver.
The data used for training the neural networks come from simulations involving the coarse solver to be corrected and a more accurate fine solver. 
We incorporated the developed deep learning framework into the parareal scheme, thereby improving its convergence property and the accuracy of the simulated wave propagation. 
Experimental results show that the presented approach dramatically enhances the stability and accuracy of the assimilated parareal scheme in a class of challenging media.   

There is much room for generalization and improvement. We presented an approach to learn the correction operator $\theta_{\Delta t^*}$ for fixed $\Delta t^*$. Would it be possible to learn $\theta_{\Delta t^*}$ as a function of $\Delta t^*$? 
The presented numerical examples comparing networks trained on $\mathtt{D_t}$ and $\mathtt{D_t^p}$ demonstrated the importance of sampling the computationally relevant, strongly causal wave fields. 
A more rigorous understanding of the wave fields manifold and devising efficient sampling algorithms would undoubtedly lead to more effective learning results.

\section*{Acknowledgment}
Nguyen and Tsai were supported partially by National Science Foundation Grant DMS-1913209. Nguyen is supported by the European Research Council Starting Grant 852821—SWING.
Tsai is also partially supported by National Science Foundation Grant DMS-1720171, and by Army Research Office, under Cooperative Agreement Number W911NF-19-2-0333. The views and conclusions contained in this document are those of the authors and should not be interpreted as representing the official policies, either expressed or implied, of the Army Research Office or the U.S. Government. The U.S. Government is authorized to reproduce and distribute reprints for Government purposes notwithstanding any copyright notation herein.

The authors thank Louis Ly for stimulating conversations on the design and training of neural networks.
The authors thank Texas Advanced Computing Center (TACC) for the computing resource. 

\bibliographystyle{abbrv}
\bibliography{main}

\appendix
\section{Appendix}

\subsection{Defining $\Lambda^\dag$}\label{Appendix:Lambda_Lambda_inv}
Given a vector $(q,p)$ and $c$, and assuming that $(q,p)=\Lambda(v,v_t)\equiv(\nabla v, v_t)$ for some function periodic $v(x,t)$. 
The main task is to find $v$. 
Of course, without any additional assumption, one can only recover a function from its partial derivatives up to a constant. So, we may define a pseudo-inverse of $\Lambda$ in the Fourier domain as:
\begin{equation} \label{eq:gradient2coordinate}
\texttt{fft}\{v\}=\begin{cases}
-i(\boldsymbol{\xi}\cdot\texttt{fft}\{ q \})|\boldsymbol{\xi}|^{-2}, & \text{for }|\boldsymbol{\xi}| \neq0,\\
C_0 & \text{for } |\boldsymbol{\xi}|=0.
\end{cases}
\end{equation}

In the setup of this paper, we need to evaluate 
$$ \Lambda^\dag\mathcal{H}_{w} \left(\Lambda \mathcal{I} \mathcal{G}_{\Delta t} \mathcal{R} [u,u_t]\right).$$
{Ideally, $C_0$ should be the integral of the wave field $u$. Since $u$ is approximated on both the fine and the coarse grid, $C_0$ has to be defined by an approximation using either grid.}
This means that $q$ comes from the first component of $\Lambda \mathcal{I} \mathcal{G}_{\Delta t} \mathcal{R} [u,u_t].$ 
We take $C_0$ to be 
More precisely, 
$$C_0 := \sum_{j=1}^{N_{\delta x}} w(x_j),$$
where $(w, w_t):=\mathcal{I}\mathcal{G}_{\Delta t} \mathcal{R} [u,u_t].$
Such choice of $C_0$ can be considered as a consistency condition, when 
$\mathcal{H}_{w}$ is replaced by the identity operator, and the coarse grid is identical to the fine grid.

\subsection{Analytic correction operator in constant media}\label{Appendix:analytic_correction_op}
Here we derive the formula for correcting the numerical dispersion error of the semi-discrete scheme in one space dimension and a constant medium. The precise context is discussed in discussed in Section~\ref{sec:correction_op}.

We use the plane wave ansatz
$$
    u(x,t) = \int A(k)e^{i(k x - \omega t)} dk + \int B(k)e^{i(k x + \omega t)} dk,
$$
where $A(k), B(k)$ are functions to be determined from the initial conditions, and $\omega=\pm ck$.

Plugging the Ansatz into the definition of energy components
\begin{eqnarray} \label{eq:energycomponent_evolution}
    \nabla u & = & \int A(k) ik e^{i(k x - \omega t)} dk + \int B(k) ik e^{i(k x + \omega t)} dk \label{eq:energycomponent_evolution1} \\
    \dfrac{\partial_t u}{c} & = & \dfrac{-i\omega}{c} \int A(k)e^{i(k x - \omega t)} dk + \dfrac{i\omega}{c} \int B(k)e^{i(k x + \omega t)} dk. \label{eq:energycomponent_evolution2}
\end{eqnarray}

Taking the Fourier transform $\mathscr{F}$ and evaluate at $t=0$, we have
\begin{eqnarray*}
    \mathscr{F}[\nabla u_0 ] &=& A(k) ik + B(k) ik \\
    \mathscr{F}[\dfrac{\partial_t u_0}{c}] & = & \dfrac{-i\omega A(k)}{c}  + \dfrac{i\omega B(k)}{c} .
\end{eqnarray*}
To keep the notation clean, let
$$
    q\equiv\mathscr{F}[\nabla u]; ~~ p\equiv\mathscr{F}[\dfrac{\partial_t u}{c}].
$$
For $k \neq 0$, rearrange the initial conditions
\begin{eqnarray*}
    \dfrac{-i q_0}{k} &=& A + B,\\
    \dfrac{i p_0 c}{\omega} &=& A - B .
\end{eqnarray*}
Solving for the coefficients, we get
\begin{eqnarray*}
    A &=& \dfrac{i}{2} \Big(-\dfrac{q_0}{k} +\dfrac{p_0 c}{\omega}\Big) \\
    B &=& \dfrac{-i}{2} \Big(\dfrac{q_0}{k} +\dfrac{p_0 c}{\omega}\Big).
\end{eqnarray*}
Then from \eqref{eq:energycomponent_evolution1},\eqref{eq:energycomponent_evolution2} we derive
\begin{equation} \label{eq:qpsolution}
    \left[\begin{array}{cc}
        q(k,t)\\
        p(k,t)
    \end{array}\right]
    = \left[\begin{array}{cc}
        \cos{\omega t} & \pm \frac{c|k|}{\omega}\sin{\omega t} \\
        \mp \frac{\omega}{c|k|} \sin{\omega t} & \cos{\omega t}
    \end{array}\right]
    \left[\begin{array}{cc}
        q_0 (k) \\
        p_0 (k)
    \end{array}\right].
\end{equation}
The $\pm$ sign reflects the two directions propagation. 
Because the wave equation has a linear dispersion relation, i.e. $\omega = \pm ck$, the matrix operator is a rotation of angle $\omega t$
\begin{equation} \label{eq:fineqpsolution}
    \left[\begin{array}{cc}
        q(k,t)\\
        p(k,t)
    \end{array}\right]
    = \left[\begin{array}{cc}
        \cos{\omega t} & \pm \sin{\omega t} \\
        \mp \sin{\omega t} & \cos{\omega t}
    \end{array}\right]
    \left[\begin{array}{cc}
        q_0 (k) \\
        p_0 (k)
    \end{array}\right].
\end{equation}
The fine solver $\mathcal{F}_{\Delta t}$ is assumed to be the exact wave propagation, so it has this linear dispersion relation. 

In contrast, we assume the coarse solver $\mathcal{G}_{\Delta t}$ solves the time-continuous space-discrete wave equation $\partial_{tt} \Tilde{u} = c^2 D^{+}_{\Delta x} D^{-}_{\Delta x} \Tilde{u}$ where $D^{\pm}_{\Delta x}$ denotes forward and backward finite difference respectively. The coarse solver introduces an error term into the above linear dispersion relation
\begin{equation}
    {\omega}_{\Delta x} = \pm ck (1-\varepsilon_{c,\Delta x}(k))
\end{equation}
for 
\begin{equation}
    \varepsilon_{c,\Delta x}(k):= ck \Big(\dfrac{k^2\Delta x^2}{24}- \dfrac{k^4 \Delta x^4}{1920} + \mathcal{O}(k\Delta x)^6 \Big).
\end{equation} 

Correspondingly, the coarse solver evolves the solution as 
\begin{equation} \label{eq:coarseqpsolution}
    \left[\begin{array}{cc}
        \Tilde{q}(k,t)\\
        \Tilde{p}(k,t)
    \end{array}\right]
    = \left[\begin{array}{cc}
        \cos{\omega_{\Delta x} t} & \pm \frac{1}{1-\varepsilon}\sin{\omega_{\Delta x} t} \\
        \mp (1-\varepsilon) \sin{\omega_{\Delta x} t} & \cos{\omega_{\Delta x} t}
    \end{array}\right]
    \left[\begin{array}{cc}
        q_0 (k) \\
        p_0 (k)
    \end{array}\right].
\end{equation}

From Equation~\eqref{eq:fineqpsolution} and Equation~\eqref{eq:coarseqpsolution}, the ideal correction operator $\Omega^*$ would bridge the gap between the coarse solver and  the fine solver. In spectral domain, the Fourier symbol of the ideal operator $\hat{\Omega}^*$ satisfies
\begin{equation} \label{eq:hatomega_approxfine}
    \hat{\Omega}^*
    \left[\begin{array}{cc}
        \cos{\omega_{\Delta x} t} & \pm \frac{1}{1-\varepsilon_{c,\Delta x}}\sin{\omega_{\Delta x} t} \\
        \mp (1-\varepsilon_{c,\Delta x}) \sin{\omega_{\Delta x} t} & \cos{\omega_{\Delta x} t}
    \end{array}\right]
    =
    \left[\begin{array}{cc}
        \cos{\omega t} & \pm \sin{\omega t} \\
        \mp \sin{\omega t} & \cos{\omega t}
    \end{array}\right].
\end{equation}

Rearranging above expression, the ideal correction operator becomes 
\begin{eqnarray*}
    \hat{\Omega}^* &= &
    \left[\begin{array}{cc}
        \cos{\omega t} & \pm \sin{\omega t} \\
        \mp \sin{\omega t} & \cos{\omega t}
    \end{array}\right]
    \left[\begin{array}{cc}
        \cos{\omega_{\Delta x} t} & \pm \frac{1}{1-\varepsilon}\sin{\omega_{\Delta x} t} \\
        \mp (1-\varepsilon) \sin{\omega_{\Delta x} t} & \cos{\omega_{\Delta x} t}
    \end{array}\right]^{-1} 
    \\
    & \simeq & 
    \left[\begin{array}{cc}
        \cos{(\omega-\omega_{\Delta x}) t} & \pm \sin{(\omega-\omega_{\Delta x}) t} \\
        \mp \sin{(\omega-\omega_{\Delta x}) t} & \cos{(\omega-\omega_{\Delta x}) t}
    \end{array}\right]\\ 
    & & + \varepsilon_{c,\Delta x} 
    \left[\begin{array}{cc}
        - \sin{\omega t} \sin{\omega_{\Delta x} t} & \mp \cos{\omega t}\sin{\omega_{\Delta x} t} \\
        \mp \cos{\omega t}\sin{\omega_{\Delta x} t} & \sin{\omega t}\sin{\omega_{\Delta x} t}
    \end{array}\right].
\end{eqnarray*}
Since $\omega - \omega_{\Delta x} = ck~\varepsilon_{c,\Delta x}(k)=\omega(c,k)\varepsilon_{c,\Delta x}(k)$.
\begin{eqnarray} 
    \hat{\Omega}^* 
    & = & \left[\begin{array}{cc}
        1 & 0 \\
        0 & 1
    \end{array}\right] + 
    \left[\begin{array}{cc}
        -\frac{(\varepsilon\omega t)^2}{2}+\mathcal{O}(\varepsilon\omega t)^4 & \pm (\varepsilon\omega t)+\mathcal{O}(\varepsilon\omega t)^3 \\
        \mp (\varepsilon\omega t)+\mathcal{O}(\varepsilon\omega t)^3  & -\frac{(\varepsilon\omega t)^2}{2}+\mathcal{O}(\varepsilon\omega t)^4
    \end{array}\right]
    +\mathcal{O}(\varepsilon), \label{eq:hatomegaapprox}
\end{eqnarray}
where $\varepsilon\omega t\equiv \varepsilon_{c,\Delta x}(k)\omega(c,k)~t$.
Applying Fourier transform to both sides of \eqref{eq:hatomega_approxfine} for physical domain, the matrix multiplications become convolutions

\begin{equation}
    \Omega^* \ast
    \left[\begin{array}{c}
        \nabla \Tilde{u} \\
        \dfrac{\partial_t \Tilde{u}}{c}
    \end{array}\right]
    = 
    \left[\begin{array}{cc}
        \mathscr{F}[\cos{\varepsilon\omega t)}] & \pm \mathscr{F}[\sin{\varepsilon\omega t}] \\
        \mp \mathscr{F}[\sin{\varepsilon\omega t}] & \mathscr{F}[\cos{\varepsilon\omega t}]
    \end{array}\right] \ast
    \left[\begin{array}{c}
        \nabla \Tilde{u} \\
        \dfrac{\partial_t \Tilde{u}}{c}
    \end{array}\right]
    \simeq
    \left[\begin{array}{cc}
        \nabla u \\
        \dfrac{\partial_t u}{c}
    \end{array}\right].
\end{equation}
{Here, the convolution and the Fourier transform are defined component-by-component.} If the term $\varepsilon \omega t$ is sufficiently small (e.g. small wave number, slow medium wave speed, and small $\Delta t$), the correction operator is a near-identity map.

\end{document}

%% file: macros.tex
\usepackage{ulem}
